\documentclass[12pt]{amsart}
\usepackage{amsmath, amsfonts, amsbsy, amsthm, amscd, graphicx}

\numberwithin{equation}{section}

\theoremstyle{plain}
\newtheorem{Theorem}{Theorem}[section]
\newtheorem{Proposition}[Theorem]{Proposition}
\newtheorem{Lemma}[Theorem]{Lemma}
\newtheorem{Corollary}[Theorem]{Corollary}

\newtheorem{Main}{Theorem}

\theoremstyle{definition}
\newtheorem{Definition}[Theorem]{Definition}

\theoremstyle{remark}
\newtheorem{Remark}{Remark}
\newtheorem{Example}{Example}
\newtheorem{Conjecture}{Conjecture}

\newcommand{\id}{\mathrm{id}}

\newcommand{\R}{{\mathbb R}}

\newcommand{\Z}{{\mathbb Z}}

\newcommand{\lk}{{\rm Lk}\,}
\newcommand{\st}{{\rm St}\,}

\newcommand{\isom}{\mathop{\rm Isom}\nolimits}

\newcommand{\cat}{\mathrm{CAT}(0)}
\renewcommand{\ker}{\mathop{\rm Ker}\nolimits}
\newcommand{\im}{\mathop{\rm Im}\nolimits}

\newcommand{\ordered}[1]{\overrightarrow{X}({#1})}
\newcommand{\orderedf}[1]{\overrightarrow{\mathcal{F}}({#1})}

\begin{document}

\title{Combinatorial Harmonic Maps and Discrete-Group Actions 
on Hadamard Spaces}

\author{Hiroyasu Izeki}
\address{Mathematical Institute, Tohoku University, 
Sendai 980-8578, Japan}
\email{izeki@math.tohoku.ac.jp}

\author{Shin Nayatani}
\address{Graduate School of Mathematics, Nagoya University,
Chikusa-ku, Nagoya 464-8602, Japan}
\email{nayatani@math.nagoya-u.ac.jp}


\dedicatory{Dedicated to Professor Takushiro Ochiai on his sixtieth
birthday}

\maketitle

\section*{Introduction.}

A harmonic map between two Riemannian manifolds is defined to be 
a critical point of the energy functional; the energy of a smooth map
is by definition
the integration of the squared norm of its differential 
over the source manifold.
By solving the heat equation associated with the energy functional,
Eells and Sampson \cite{EelSam} proved that if two manifolds are compact
and the target has nonpositive sectional curvature, any smooth map 
can be deformed to a harmonic map.
This existence result was later extended to the equivariant setting 
\cite{Cor1, Don, JosYau1, Lab}.

Application of the harmonic map theory to the superrigidity of lattices
started with the work of Corlette \cite{Cor2}, in which it was proved that 
superrigidity over archimedian fields holds for lattices in $Sp(n,1)$ 
($n\geq 2$) and $F_4^{-20}$.
By developing the relevant harmonic map theory, Gromov and Schoen 
\cite{GroSch} carried the preceding result over to the $p$-adic
case, and thereby established the arithmeticity of such lattices.
These works were followed by those of Mok-Siu-Yeung \cite{MSY} and Jost-Yau 
\cite{JosYau2}, which established a result including Corlette's one
mentioned 
above as well as the cocompact case of 
the superrigidity theorem of Margulis \cite{Mar1} concerning lattices 
in real semisimple Lie groups of rank $\geq 2$.

Wang \cite{Wan1, Wan2} has taken an important step in developing 
an analogy of the story 
described hitherto when the source and target spaces are respectively 
a simplicial complex and a Hadamard space, that is, a complete CAT($0$) 
space; this is the setting necessary for the geometric-variational approach 
to the Margulis superrigidity for lattices in semisimple algebraic groups
over $p$-adic fields.
Wang formulated the notion of energy
of an equivariant map between the spaces of the above sort, 
and established the existence of an energy-minimizing equivariant map
assuming that the target space is locally compact and the action is reductive.
He then gave an application to the isometric action of a discrete group,
arising as the covering transformation group of a finite simplicial complex,
on a Hadamard space.
Under the same assumptions as above, he formulated a general criterion 
for the existence of a fixed point 
in terms of a Poincar\'e-type constant for maps from the link
of a vertex of the source to a tangent cone of the target.
In this paper we shall refer to this constant as {\em Wang's invariant}.
He also proved a more concrete fixed-point theorem when the target space 
was a Hadamard manifold.
It should be mentioned that there are also various forms of harmonic map 
theories with singular source spaces \cite{Jos1, Jos2, KotSun, Leb}.

The purpose of the present paper is to push forward with Wang's work,
especially when the target Hadamard space is {\em non-locally compact}
and/or {\em singular}.
We shall thereby establish further applications of combinatorial harmonic
map theory to the isometric actions of discrete groups on Hadamard spaces.
Our main result strengthens Wang's result mentioned above, to the effect that 
we drop both of the 
assumptions that the target space is locally compact and that the action 
is reductive.

\begin{Main}\label{main_result1}
Let $X$ be a connected simplicial complex. 
We assume that the links of all vertices of $X$ are connected.
Let $\Gamma$ be a finitely generated group acting by automorphisms, 
properly discontinuously and cofinitely on $X$. 
Let $Y$ be a Hadamard space.
Suppose that Wang's invariant is not less than a constant $C$ 
$> 1/2$ for all vertices of $X$ and all points of $Y$.
Then any isometric action of $\Gamma$ on $Y$ has a fixed point.
\end{Main}

By estimating Wang's invariant in some instances, we deduce the following 
concrete fixed-point theorems from Theorem \ref{main_result1}.

\begin{Main}\label{main_result2}
Let $X$ and $\Gamma$ be as in Theorem \ref{main_result1}.
We assume that the first nonzero eigenvalues of the Laplacians for
the links of all vertices of $X$ are greater than $1/2$.
Let $Y$ be one of the following spaces{\rm :}
\begin{enumerate}
\renewcommand{\theenumi}{\roman{enumi}}
\renewcommand{\labelenumi}{(\theenumi)}
\item a Hadamard space all of whose tangent cones are isometric 
to closed convex cones of Hilbert spaces,
\item an $\R$-tree,
\item a product of Hadamard spaces which are either of type {\rm (i)} 
or {\rm (ii)}.
\end{enumerate}
Then any isometric action of $\Gamma$ on $Y$ has a fixed point.
\end{Main}

\begin{Main}\label{main_result3}
Let $X$ and $\Gamma$ be as in Theorem \ref{main_result1}.
Let $Y$ be the Euclidean building 
$PGL(3,\mathbb{Q}_2)/PGL(3,\mathbb{Z}_2)$.
Then there is a constant $C\in(1/2,1)$ such that if the first nonzero 
eigenvalues of the Laplacians for the links of all vertices of $X$ are 
greater than $C$,
then any isometric action of $\Gamma$ on $Y$ has a fixed point.
\end{Main}

While completing the earlier version of this paper, we learned that 
Gromov \cite{Gro2} also described similar approach to fixed-point 
theorems.
Here we summarize the results in \cite{Gro2} which are most closely 
related to those in this paper.
He proves the assertion of Theorem \ref{main_result1} by using 
the scaling limit argument and assuming that the condition 
on Wang's invariant as in Theorem \ref{main_result1} (or its variant)
holds for various scaling Hausdorff limits of the target 
Hadamard space \cite[p.~128]{Gro2}.
As a corollary, he formulates a fixed-point theorem for the class 
$\overline{\mathcal{C}}_{\rm reg}$ of Hadamard spaces, that is 
the minimal class of Hadamard spaces containing all Hadamard manifolds
and that is closed under taking totally geodesic subspaces
and under Hausdorff limits.
He gives a few more results on the existence of a fixed-point 
\cite[p.~118, p.~125]{Gro2} under assumptions which assure the exponential 
decay of the energy along his discrete heat flow and thereby conclude the 
convergence of the flow to a constant map.

Our proof of Theorem \ref{main_result1} is 
along more direct approach without referring to scaling Hausdorff limits
of the target, thereby avoiding the additional assumption put in Gromov's result.
It should also be mentioned that the Hadamard spaces in Theorem \ref{main_result2}
are included in Gromov's class $\overline{\mathcal{C}}_{\rm reg}$.
On the other hand, as pointed out by Gromov, most Hadamard spaces 
are not in his class $\overline{\mathcal{C}}_{\rm reg}$; in particular,
this is the case for most Euclidean buildings including the one in 
Theorem \ref{main_result3}.

We shall now give an overview of this paper.
In \S 1 we collect some definitions and results on Hadamard spaces.
In \S 2 we define the energy of an equivariant map from a simplicial
complex into a Hadamard space.
Here an action of a discrete group on the source and a homomorphism 
of it into the isometry group of the target are given, and the equivariance
is with respect to this homomorphism.
We derive the first variation formula for the energy
(Proposition  \ref{variation_formula}).
It follows that an energy-minimizing equivariant map, if exists,
satisfies a certain balancing condition, which will be precisely
stated using the notion of barycenter on the tangent cones of the
target.
We adopt this last condition as the defining property of a harmonic 
equivariant map.
We then state an existence result for an energy-minimizing equivariant 
map (Theorem \ref{existence}).
It is worth while to emphasize that our formulation is more general than
Wang's, to the effect that we do not require the action of the discrete group 
on the source is free.
Our results thus apply to wider range of examples of discrete groups.
We also mention that Gromov gives an argument which reduces the case of 
non-free action to that of free one
by enlarging the source simplicial complex.
In \S 3 we shall write down Bochner-type formulas for an equivariant
map (Proposition \ref{nonlinear_case}).
Such a formula is necessary in deriving an estimate of the gradient
of the energy functional in \S 4.
It also explains the way that Wang's criterion forces an energy-minimizing
equivariant map to be a constant map.

In \S 4 we present a different scheme to find a constant equivariant map 
and thus prove the existence of a fixed point for the action.
Under the assumption on Wang's invariant as in Theorem \ref{main_result1},
we prove that any equivariant map is deformed to a constant map 
along the gradient flow associated with our energy functional.
We therefore prove Theorem \ref{main_result1}. 
It should be mentioned that the gradient flow we employ here 
was introduced by Jost \cite{Jos4} and Mayer \cite{May} in a more general setting, 
and we shall use their results in the proof.

In \S 5 we shall give lower and upper estimates of Wang's invariant
in terms of the first nonzero eigenvalue of the Laplacian for the vertex' link
of the source and a certain geometric invariant for the tangent cone of 
the target (Definition \ref{delta}).
The former invariant is rather well studied partly because of its
relation to Kazhdan's property (T).
On the other hand, the latter invariant, which we call $\delta$ here, 
is newly introduced in this paper as far as the authors know.
It measures the extent to which a tangent cone of a Hadamard space
differs from a Euclidean space, regarding the property of their barycenter 
correspondences.
This invariant looks measuring some sort of rigidity of a tangent cone
of a Hadamard space, and might be useful in other contexts studying
such spaces.
It takes values in the interval $[0,1]$, and if it is zero, then Wang's 
invariant coincides 
with the first nonzero eigenvalue of the Laplacian for the vertex' link, 
which is computable in principle.
We thereby generalize Wang's fixed-point theorem for a Hadamard manifold
to one for Hadamard spaces all of whose tangent cones have vanishing $\delta$
(Theorem \ref{flexible_case}).
Other than the trivial examples of Hadamard manifolds, we prove that
$\mathbb{R}$-trees also satisfy this condition and also that products of 
known examples give rise to new ones.
In particular, we obtain Theorem \ref{main_result2}.
The corollary has an interesting consequence to the effect that there are 
many discrete groups which admit no properly discontinuous isometric 
action on any Hadamard manifold.
This provides a negative answer to the question asked by Gromov 
in \cite{Gro1}.
There are certainly many examples of Hadamard spaces for which $\delta$
is positive for some of their tangent cones.
In fact, this is the case for most Euclidean buildings.
Our criterion is particularly meaningful if $\delta$ is small for all 
tangent cones.
Though we cannot compute $\delta$ for Euclidean buildings at this time, 
in \S 6 we give a reasonably good estimate of $\delta$ for a certain 
Euclidean building. 
Theorem \ref{main_result3} is an immediate consequence of this estimate
and Theorem \ref{main_result1}.

In Appendix we write down a Bochner-Weitzenb\"ock-type formula for
an equivariant cochain of an arbitrary degree on a simplicial complex
(Proposition \ref{linear_case}).
The formula for a cochain of degree one is referred in the proof of 
Proposition \ref{nonlinear_case}.
The formula is essentially contained in the work of Ballmann
and \'Swi\c{a}tkowski \cite{BalSwi}, where they used it to prove 
that certain groups have Kazhdan's property (T).
We expect, however, that our presentation is more transparent and direct.
It should be mentioned that the formula of this type has its origin in the 
work of Garland \cite{Gar} where the case of Euclidean buildings is treated.
Pansu \cite{Pan} gave a geometric explanation to Garland's computation
when the building is of type $\widetilde{A}_2$.

Theorems 1, 2, 3 will be restated as Theorem
\ref{fixed_point_theorem1}, Corollary \ref{fixed_point_theorem2}
and Corollary \ref{fixed_point_theorem_for_A_2_tilde_building}
respectively.

Part of this work was announced in \cite{IzeNay1, IzeNay2}.

\smallskip\noindent
{\bf Acknowledgements.}\quad
We would like to thank G\'erard Besson, Koji Fujiwara, Masahiko Kanai,
Shin Kato, Koichi Nagano, Shin-ichi Ohta and Andrzej \.Zuk for helpful 
discussions, suggestion and interest in this work.

\section{Preliminaries on Hadamard spaces}

In this section, we collect some definitions and results concerning
$\cat$ and Hadamard spaces.
Most of them are well-known and available in literature; we refer 
the reader to \cite{BriHae} for a detailed exposition on these spaces.
Possible exceptions are Lemmas \ref{angle2} and \ref{fixed-point_set},
to which we give complete proofs.

Let $Y$ be a metric space and $p,q \in Y$.
A {\it geodesic} joining $p$ to $q$ is 
an isometric embedding $c$ of a closed interval $[0,l]$ into
$Y$ such that $c(0)=p$ and $c(l)=q$. 
We say that $Y$ is a {\it geodesic space} if any two points in $Y$ 
are joined by a geodesic.  

Consider a triangle in $Y$ with vertices $p_1,p_2,p_3 \in Y$ and 
three geodesic segments $p_1p_2, p_2p_3, p_3p_1$ joining them. 
We denote this triangle by $\Delta(p_1,p_2,p_3)$ and
call such a triangle a {\it geodesic triangle}.  
Take
a triangle $\Delta(\overline{p_1},\overline{p_2},\overline{p_3})$ in
$\R^2$ so that $d_Y(p_i, p_j)=d_{\R^2}(\overline{p_i},\overline{p_j})$.
We call $\Delta(\overline{p_1},\overline{p_2},\overline{p_3})$ a {\it
comparison triangle} for $\Delta(p_1,p_2,p_3)$. 
A point $\overline{q}\in \overline{p_i}\overline{p_j}$
is called a {\it comparison point} for $q\in p_ip_j$ if 
$d_Y(p_i,q) = d_{\R^2}(\overline{p_i}, \overline{q})$.
If $d_Y(q_1,q_2) \leq d_{\R^2}(\overline{q_1}, \overline{q_2})$ for any
pair of points $q_1, q_2$ on the sides of $\Delta(p_1,p_2,p_3)$ and
their comparison points $\overline{q_1}, \overline{q_2}$, 
then $\Delta(p_1,p_2,p_3)$ is said to satisfy the {\it $\cat$ inequality}. 
If every geodesic triangle in $Y$ satisfies the $\cat$ inequality, 
then $Y$ is called a {\it $\cat$ space}. 

\vskip0.5cm
\hskip2.5cm
\includegraphics{cat0.1}
\vskip0.5cm

\noindent
If $Y$ is a $\cat$ space, it
is easy to verify that $Y$ has two important properties
which we shall use frequently in what follows: the uniqueness of 
a geodesic connecting given two points and the contractibility. 
If a $\cat$ space
$Y$ is complete as a metric space, it is called a {\it Hadamard
space}. 

\begin{Proposition}
\label{barycenter}
Let $Y$ be a Hadamard space, 
let $p_1,\dots,p_m$ be points of $Y$,
and let $t_1,\dots,t_m$ be real numbers satisfying $t_i\geq 0$
and $\sum_{i=1}^m t_i = 1$.
Then there exists a unique point $p\in Y$
which minimizes the function
$$
F(q) = \sum_{i=1}^m t_i d(p_i,q)^2,\quad q\in Y.
$$
\end{Proposition}
For the proof, see \cite[p. 639, Lemma 2.5.1] {KorSch1}. 
We call the point $p$ the \textit{barycenter of} 
$\{p_1, \dots ,p_m\}$ \textit{with weight} $\{t_1, \dots, t_m\}$.  
If $t_i = 1/m$ for all $i$, we simply call it the \textit{barycenter of}   
$\{p_1, \dots, p_m\}$. 

\begin{Definition}\label{tangent cone}
Let $Y$ be a Hadamard space.\\
(1)\,\, Let $c$ and $c'$ be two nontrivial geodesics in $Y$ starting 
from $p \in Y$. The {\it angle} $\angle_p(c,c')$ between $c$ and $c'$ 
is defined by 
\begin{equation*}
 \angle_p(c,c')= \lim_{t,t' \rightarrow 0}
 \angle_{\overline{p}}(\overline{c(t)},\overline{c'(t')}),
\end{equation*}
 where
 $\angle_{\overline{p}}(\overline{c(t)},\overline{c'(t')})$ 
 denotes the angle between the sides
 $\overline{p}\overline{c(t)}$ and $\overline{p}\overline{c'(t)}$
 of the comparison triangle
 $\Delta(\overline{p},\overline{c(t)},\overline{c'(t')})
 \subset \R^2$.\\
(2)\,\, Let $p \in Y$.
 We define an equivalence relation  $\sim$ on the set of
 nontrivial geodesics starting from $p$ by 
 $c \sim c' \Longleftrightarrow \angle_p(c,c')=0$. 
 Then the angle $\angle_p$ induces a distance on the quotient 
 $(S_pY)^{\circ} = \{\text{nontrivial geodesics starting from } 
 p \}/\sim$, 
 which we denote by the same symbol $\angle_p$.  The completion 
 $(S_pY,\angle_p)$ of the metric space $((S_pY)^{\circ}, \angle_p)$ is
 called the {\it space of directions} 
 at $p$.\\  
(3)\,\, Let $TC_pY$ be the cone over $S_pY$, namely,
\begin{equation*}
 TC_pY = (S_pY \times \R_+) / (S_pY \times \{0\}). 
\end{equation*}
 Let $W,W' \in TC_pY$. We may write $W=(V,t)$ and $W'=(V',t')$, where
 $V,V' \in S_pY$ and $t,t' \in \R_{+}$.   Then
\begin{equation*}
 d_{TC_pY}(W, W')= t^2 + {t'}^2 - 2tt'\cos \angle_p(V,V')
\end{equation*}
 defines a distance on $TC_pY$. The metric space $(TC_pY, d_{TC_pY})$ is
 again a Hadamard space and is called the {\it tangent cone} at $p$.  We
 define an ``inner product''  on  $TC_pY$ by
\begin{equation*}
 \langle W, W' \rangle = tt'\cos \angle_p(V,V').
\end{equation*}
 We shall often denote the length $t$ of $W$ by $|W|$; thus we have
 $|W|=\sqrt{\langle W,W \rangle}=d_{TC_pY}(0_p,W)$, where $0_p$ denotes
 the origin of $TC_pY$.   \\
(4)\,\, Define a map $\pi_p : Y \longrightarrow TC_pY$ by
 $\pi_p(q)=([c], d_Y(p,q))$, where $c$ is the geodesic joining $p$ to
 $q$ and $[c]\in S_pY$ is the equivalence class of $c$. Then $\pi_p$ is
 distance non-increasing.  
\end{Definition}

A complete, simply connected Riemannian manifold $Y$ with nonpositive
sectional curvature, called a {\it Hadamard manifold}, is a typical 
example of Hadamard space. 
For such a $Y$, $S_pY$ (resp. $TC_pY$) is the unit tangent sphere
(resp. the tangent space) at $p$. 
The map $\pi_p$ is the inverse of the exponential map, and it is 
well-known that the exponential map is distance non-decreasing 
for such a $Y$.

The following results will be used in the proof of 
Proposition \ref{variation_formula}.

\begin{Lemma}{\cite[Corollary II 3.6]{BriHae}}\label{angle}
 Let $Y$ be a Hadamard
 space and $p \in Y$. Take two nontrivial geodesics $c$, $c'$ starting
 from $p$, and fix a point $q$ on $c'$.  Then we have
\begin{equation*}
  \cos \angle_p (c,c') = \lim_{t \rightarrow 0} 
 \frac{d_Y(p,q) - d_Y(c(t),q)}{t}.
\end{equation*}
\end{Lemma}

\begin{Lemma}\label{angle2}
 Let $Y$ be a Hadamard space and $p,q \in Y$.  Let $c(t)$
 $($resp.~$c'(t)$\,$)$ be a nontrivial geodesic starting from $p$
 $($resp.~$q$\,$)$. Then we have
 \begin{equation*}
  \cos \angle_p(c,pq) = \lim_{t\rightarrow 0} 
       \frac{ d_Y(p,c'(t))-d_Y(c(t),c'(t))}{t}, 
 \end{equation*}
 where $pq$ denotes the geodesic starting from $p$ and terminating at
 $q$. 
\end{Lemma}

\begin{proof}
 Consider the subembedding 
 $(\overline{q},\overline{p},\overline{c(t)},
 \overline{c'(t)}) \subset \R^2$ 
 of $(q,p,c(t),c'(t))$ (see \cite[II 1.10]{BriHae}) such that
 \begin{equation*}
  \begin{split}
   & d_Y(p,q)=d_{\R^2}(\overline{p},\overline{q}), \quad
   d_Y(p,c(t))=d_{\R^2}(\overline{p},\overline{c(t)}), \\
   & d_Y(c(t),c'(t))=d_{\R^2}(\overline{c(t)},\overline{c'(t)}), \quad
   d_Y(c'(t),q)=d_{\R^2}(\overline{c'(t)},\overline{q}). 
  \end{split}
 \end{equation*}
 Then we have 
 \begin{equation*}
  d_Y(p,c'(t))\leq d_{\R^2}(\overline{p},\overline{c'(t)}), \quad
  d_Y(q,c(t))\leq d_{\R^2}(\overline{q},\overline{c(t)}). 
 \end{equation*}
 Set 
 \begin{equation*}
   a(t)=d_{\R^2}(\overline{c(t)},\overline{c'(t)}), \quad
   b(t)=d_{\R^2}(\overline{q},\overline{c(t)}), \quad
   \widetilde{a}(t)=d_{\R^2}(\overline{p},\overline{c'(t)}). 
 \end{equation*} 
 Note that
 $b(0)=d_{\R^2}(\overline{p},\overline{q})=d_Y(p,q)$. 
 By the cosine rule, we have
  \begin{eqnarray*}
  \widetilde{a}(t) &=& b(0)   - t\cos
  \angle_{\overline{q}}(\overline{q}\overline{p}, 
          \overline{q}\overline{c''(t)})  + o(t), \\
  b(t) &=& b(0)   - t\cos
  \angle_{\overline{p}}(\overline{p}\overline{q}, 
          \overline{p}\overline{c(t)})  + o(t). 
 \end{eqnarray*}
The latter expression gives 
 \begin{eqnarray*}
   a(t)^2 & = & t^2 + b(t)^2 - 2tb(t)
   \cos \angle_{\overline{q}}
        (\overline{q}\overline{c(t)}, \overline{q}\overline{c'(t)}) \\
   &=& b(0)^2 - 2tb(0)\left(\cos \angle_{\overline{p}}
        (\overline{p}\overline{q}, \overline{p}\overline{c(t)}) + 
         \cos \angle_{\overline{q}}
        (\overline{q}\overline{p}, \overline{q}\overline{c'(t)})\right)
      + o(t),
 \end{eqnarray*}
 and therefore
 \begin{equation*}
  a(t) = b(0) - t\left(\cos \angle_{\overline{p}}
        (\overline{p}\overline{q}, \overline{p}\overline{c(t)}) + 
         \cos \angle_{\overline{q}}
        (\overline{q}\overline{p}, \overline{q}\overline{c'(t)})\right)
        + o(t).  
 \end{equation*}
 Hence we get
 \begin{equation*}
  a(t)+b(0)-\widetilde{a}(t)-b(t)= o(t). 
 \end{equation*}
  By Lemma \ref{angle}, 
 \begin{eqnarray*}
  \cos \angle_p(c,pq) & = 
  &  \lim_{t \rightarrow 0}\frac{d_Y(p,q)-d_Y(c(t),q)}{t} \geq 
     \lim_{t \rightarrow 0}\frac{b(0)-b(t)}{t} \\
  & = & \lim_{t \rightarrow 0}\frac{\widetilde{a}(t)-a(t)}{t} \geq
 \lim_{t \rightarrow 0}\frac{d_Y(p,c'(t))-d(c(t),c'(t))}{t}.
 \end{eqnarray*}
 
 Take another subembedding
 $(\overline{q}',\overline{c(t)}',\overline{p}',\overline{c'(t)}') 
 \subset \R^2$ 
 of $(q,c(t),p,c'(t))$ such that
 \begin{equation*}
  \begin{split}
   & d_Y(p,c(t))=d_{\R^2}(\overline{p}',\overline{c(t)}'), \quad
   d_Y(c'(t),q)=d_{\R^2}(\overline{c'(t)}',\overline{q}'), \\  
   & d_Y(p,c'(t))= d_{\R^2}(\overline{p}',\overline{c'(t)}'), \quad
   d_Y(q,c(t))= d_{\R^2}(\overline{q}',\overline{c(t)}'). 
  \end{split}
 \end{equation*}
 Then we have  
 \begin{equation*} 
 d_Y(p,q)\leq d_{\R^2}(\overline{p}',\overline{q}'), \quad
  d_Y(c(t),c'(t))\leq d_{\R^2}(\overline{c(t)}',\overline{c'(t)}').   
 \end{equation*}
 By the same computation as above, we get the opposite
 inequality 
  \begin{equation*}
  \cos \angle_p(c,pq) \leq \lim_{t\rightarrow 0} 
       \frac{ d_Y(p,c'(t))-d_Y(c(t),c'(t))}{t}.
 \end{equation*}
 This completes the proof. 
\end{proof}

We denote the isometry group of $Y$ by $\isom(Y)$.
Let $G$ be a subgroup of $\isom(Y)$ which fixes a point $p \in Y$.  
Note that the action of $G$ on $Y$ induces an isometric action on $S_pY$. 
And this action extends to one on $TC_pY$ in a natural way.  

\begin{Lemma}\label{fixed-point_set}
 Let $G$ be a finite subgroup of $\isom (Y)$ which fixes a
 point $p \in Y$. 
 Denote by $Y^G$ $($resp.~$(TC_pY)^G$\,$)$ the fixed-point set
 with respect to the action of $G$ on $Y$ $($resp.~$TC_pY$\,$)$. Then 
 $(TC_pY)^G$ coincides with the closure of 
 $\mathbb{R}_+ \pi_p(Y^G) = \{ tV \mid t\geq 0, V\in \pi_p(Y^G) \}$.
\end{Lemma}

\begin{proof} It is clear that $(TC_pY)^G \supset \mathbb{R}_+ \pi_p(Y^G)$, 
and that
 $(TC_pY)^G$ is a closed subset of $TC_pY$.  Thus it suffices to show
 that $(TC_pY)^G$ is contained in the closure of $\mathbb{R}_+ \pi_p(Y^G)$.  
 Let $W \in (TC_pY)^G$.  
We shall find a sequence of points $\{ q_j \}_{j=0}^{\infty} \subset Y^G$,
$q_j \neq p$, 
which satisfy $\angle_p([c_j],W)\rightarrow 0$ if $c_j$ is the geodesic
joining $p$ to $q_j$.
Then the sequence $\{W_j=([c_j], |W|)\}_{j=0}^{\infty}$ 
converges to $W$, and hence $W$ belongs to the closure of 
$\mathbb{R}_+ \pi_p(Y^G)$. 

 First note that there is a sequence of geodesics 
 $\{c'_j\}_{j=0}^{\infty}$ such that 
 $\angle_p([c'_j],W) \rightarrow 0$. For such a sequence $\{c'_j\}$, we
 must have  
 $\max_{\gamma \in G}\angle_p(c'_j, \gamma c'_j) \rightarrow 0$. 
 Define $R_{j,t}$ to be the radius of $Gc'_j(t)$, that is,
 \begin{equation*}
  R_{j,t} = \inf \{r>0 \mid Gc'_j(t) \subset B(q,r) \text{ for some }q \in
  Y\}, 
 \end{equation*}
 where $B(q,r)$ denotes the metric ball with center $q$ and 
 radius $r$. Then there exists a unique point $q_{j,t}$ called the {\it
 circumcenter} of $Gc'_j(t)$ that satisfies 
 $Gc'_j(t)\subset B(q_{j,t},R_{j,t})$ 
 (\cite[Proposition II 2.7]{BriHae}).
Note that $q_{j,t}\in Y^G$.
 It is clear that 
 $$ \max_{\gamma \in G}d_Y(c'_j(t),\gamma c'_j(t))\geq R_{j,t} $$
 holds, and we have 
 \begin{equation*}
 \max_{\gamma \in G} \angle_{\overline{p}}(\overline{p}
 \overline{c'_j(t)},
 \overline{p}\overline{\gamma c'_j(t)}) \geq 2 \sin^{-1}
 \frac{R_{j,t}}{2t}, 
 \end{equation*}
 where the angles in the left-hand side are the comparison angles.  Set
 $\varepsilon_{j} = \max_{\gamma \in G}\angle_p(c'_j,\gamma c'_j)$. Then 
 $\max_{\gamma \in G}\angle_{\overline{p}}(\overline{p}\overline{c'_j(t)},
 \overline{p}\overline{\gamma c'_j(t)}) \rightarrow
 \varepsilon_j$
 as $t \rightarrow 0$, 
 and $\varepsilon_j \rightarrow 0$ as $j \rightarrow \infty$ by our
 choice of $c'_j$. 
 Therefore we can take a sequence $\{t_j\}$ so that $t_j \searrow 0$ and  
 $\max_{\gamma \in G}
 \angle_{\overline{p}}(\overline{p}\overline{c'_j(t_j)}, 
 \overline{p}\overline{\gamma c'_j(t_j)})\rightarrow 0$
 as $j\rightarrow \infty$. 
 Let $R_j=R_{j,t_j}$ and $q_j=q_{j,t_j} \in Y^G$.  Then we have 
 $\sin^{-1} R_j/2t_j \rightarrow 0$ as
 $j \rightarrow \infty$.  In particular, $R_j/t_j \rightarrow 0$, and 
therefore we may assume $q_j \neq p$ for each $j$.  
 On the other hand, we have
 \begin{equation}
 \label{angle_goes_to_zero}
 \angle_{\overline{p}}(\overline{p}\overline{c'_j(t_j)},
 \overline{p}\overline{q_j}) \leq \sin^{-1} \frac{R_j}{t_j},
 \end{equation}
 where the right-hand side is the Euclidean angle 
 $\angle_{a} (ab,ac)$ at the vertex $a$ of the triangle
 $\Delta(a,b,c)$ with 
 $d_{\R^2}(a,b)=t_j$, $d_{\R^2}(b,c)=R_j$ and
 $\angle_{c}(ca,cb)= \pi/2$. 
Denote the geodesic joining $p$ to $q_j$ by $c_j$.  
Then by \eqref{angle_goes_to_zero}, we see 
that $\angle_p([c'_j],[c_j]) \rightarrow 0$ as $j \rightarrow \infty$.
Therefore $\angle_p([c_j],W)\rightarrow 0$.  
This completes the proof. 
\end{proof}

\section{Energy of equivariant maps}

In this section, after some preliminaries on simplicial complexes,
we shall define the energy of an equivariant map from a simplicial
complex to a Hadamard manifold.
We compute the first variation of the energy, which gives a
necessary condition that an energy-minimizing equivariant map must satisfy,
and also motivates the definition of harmonicity of an equivariant map.
We shall conclude with an existence result for an energy-minimizing 
equivariant map.

Let $X$ be a simplicial complex.
Throughout this paper, we shall assume that $X$ is connected, locally finite 
and finite-dimensional.
For each $k\geq 0$, let $X(k)$ denote the set of 
$k$-simplices in $X$.  
For a $k$-simplex $s \in X(k)$ and $l>k$, let $X(l)_s$ denote the set
of $l$-simplices containing $s$.  
It is often convenient to deal with ordered simplices instead of 
simplices without ordering.  We denote the set of ordered $k$-simplices
by $\overrightarrow{X}(k)$, and an element of $\overrightarrow{X}(k)$ by  
$(x_0,x_1,\dots,x_k)$. For $s \in \ordered{k}$ and $l>k$, we set
\begin{equation*}
 \ordered{l}_s = \{t=(x_0,\dots,x_k,\dots x_l) \in \ordered{l}
                    \mid (x_0,\dots, x_k)=s \}. 
\end{equation*}

\begin{Definition} 
\label{admissible_weight}
 A positive function 
 $m:\bigcup_{k \geq 0} X(k) \longrightarrow \R$ is called an 
 {\it admissible weight} if it satisfies
\begin{equation}
\label{weight}
 m(s) = \sum_{t \in X(k+1)_s} m(t)
\end{equation}
for all $s \in X(k)$ such that $X(k+1)_s \neq \emptyset$.
We often regard $m$ as a symmetric function,
still denoted by the same symbol $m$, on the set of ordered simplices 
$\bigcup_{k \geq 0} \ordered{k}$. 
\end{Definition}

Given a positive function on the set of maximal simplices of $X$, 
we can define an admissible weight inductively by \eqref{weight}. 
We call the admissible weight $m$ defined by setting $m(s)=1$ for every
maximal simplex $s$ of $X$ the {\it standard admissible weight} of $X$.  
Throughout this paper, we shall assume that $X$ is equipped with an 
admissible weight $m$.
By an {\it automorphism} of $X$, we mean a simplicial automorphism of $X$ that preserves $m$. 
Note that this last requirement is automatic if $m$ is the standard 
admissible weight.

Let $\Gamma$ be a finitely generated (hence, countable) group acting 
by automorphisms and properly discontinuously on $X$. 
Since the action of $\Gamma$ is simplicial, $\Gamma$ acts on $X(k)$ and 
$\overrightarrow{X}(k)$ for each $k$. 
We denote by $\mathcal{F} (k) \subset X(k)$ (resp. 
$\overrightarrow{\mathcal{F}}(k)\subset \overrightarrow{X}(k))$ 
a set of representatives of 
the $\Gamma$-orbits in $X(k)$ (resp. $\overrightarrow{X}(k)$). 
The isotropy subgroup of an ordered or unordered simplex $s$ is
denoted by $\Gamma_s$.  
For $s \in \ordered{k}$ and $l>k$, we denote a set of representatives of 
the $\Gamma_s$-orbits in $\ordered{l}_s$ by $\orderedf{l}_s$.    

The action of $\Gamma$ is called \textit{cofinite} if $\mathcal{F} (k)$
is finite for each $k$.  
Throughout this paper, we shall assume, for the sake of simplicity, 
that the action of $\Gamma$ is cofinite. 

Recall that the \textit{star} $\st s$ of $s\in X(k)$ is the subcomplex 
of $X$ such that $t \in \st s$ if and only if $s \cup t$ is a simplex of 
$X$.  
The \textit{link} of $s\in X(k)$, denoted by $\lk s$, is the
subcomplex of $\st s$ such that $t \in \lk s$ if and
only if $t \in \st s$ and $t$ is disjoint from $s$. 
The action of $\Gamma$ on $X$ is called \textit{very free} if 
$\st x \cap \gamma \st x = \emptyset$ for any $\gamma \in \Gamma
\setminus \{\id\}$ and $x \in X(0)$.
It is easy to check that
if the action of $\Gamma$ is very free, then the quotient 
$\Gamma \backslash X$ inherits a structure of simplicial complex from
$X$.   

\begin{Lemma}
\label{sum}
 Let $\psi$ be a $\Gamma$-invariant function on $\ordered{l}$, and 
 $0 \leq k \leq l-1$. Then we have
 \begin{equation*}
  \sum_{u \in \orderedf{l}}\frac{1}{|\Gamma_u|}\psi(u)
  =    \sum_{s \in \orderedf{k}} \frac{1}{|\Gamma_s|}  
       \sum_{u \in \ordered{l}_s} \psi(u).
 \end{equation*}
\end{Lemma}

\begin{proof} First note that the both sides of the formula
  are independent of the choice of the representative sets
 $\orderedf{k}$, $\orderedf{l}$. 
It is easy to verify that, starting with a given $\mathcal{F}(0)$, 
 we can choose representative sets inductively so that
\begin{equation*}
\orderedf{k} = \bigcup_{t \in \orderedf{k-1}} \orderedf{k}_t, \quad
 k\geq 1. 
\end{equation*}
On the other hand, for $t \in \ordered{l-1}$ and $u \in \ordered{l}_t$, 
we have 
$\# \{\gamma u \ |\ \gamma \in \Gamma_t \}= 
\left.|\Gamma_t|\right/|\Gamma_u|$.  
Thus we may rewrite
\begin{eqnarray*}
       \sum_{u \in \orderedf{l}}\frac{1}{|\Gamma_u|}\psi(u)
  &=&  \sum_{t \in \orderedf{l-1}} 
       \sum_{u \in \orderedf{l}_t}\frac{1}{|\Gamma_u|}\psi(u) \\
  &=&  \sum_{t \in \orderedf{l-1}}
       \sum_{u \in \ordered{l}_t} \frac{|\Gamma_u|}{|\Gamma_t|}
       \frac{1}{|\Gamma_u|}\psi(u) \\
  &=&  \sum_{t \in \orderedf{l-1}}
       \frac{1}{|\Gamma_t|} \sum_{u \in \ordered{l}_t}
       \psi(u) \\
  &=&  \sum_{s \in \orderedf{l-2}}
       \sum_{t \in \ordered{l-1}_s}
       \frac{|\Gamma_t|}{|\Gamma_s|}\frac{1}{|\Gamma_t|}
       \sum_{u \in \ordered{l}_t} \psi(u) \\
  &=&  \sum_{s \in \orderedf{l-2}} \frac{1}{|\Gamma_s|}  
       \sum_{u \in \ordered{l}_s} \psi(u),
\end{eqnarray*}
and so on.
This completes the proof. 
\end{proof}

Suppose that a homomorphism 
$\rho:\Gamma \longrightarrow \isom(Y)$ is given, where $Y$ is a Hadamard space.
A map $f: X(0)\longrightarrow Y$ is said to be
$\rho$-\textit{equivariant} if $f$ satisfies 
$f(\gamma x) = \rho(\gamma)f(x)$ for every $\gamma \in \Gamma$ and $x \in X(0)$.  
Now we define the energy of 
a $\rho$-equivariant map, which is a nonlinear analogue of the $L^2$-norm
of $df$ for an equivariant $0$-cochain $f$ (see Appendix). 

\begin{Definition}
 \label{energy}
 Let $f:X(0)\longrightarrow Y$ be a $\rho$-equivariant map, and $x \in X(0)$.
 We define the \textit{energy density} $e_{\rho}(f)(x)$ by 
  \begin{equation*}
     e_{\rho}(f)(x) = \frac{1}{2}
      \sum_{(x,y)\in \overrightarrow{\mathcal{F}}(1)_x} 
       \frac{m(x,y)}{|\Gamma_{(x,y)}|}   d_Y(f(x), f(y))^2.
  \end{equation*}
 Then the \textit{energy} $E_{\rho}(f)$ of $f$ is defined 
 by
  \begin{equation*}
     E_{\rho}(f) = \sum_{x \in \mathcal{F}(0)} e_{\rho}(f)(x). 
  \end{equation*}
\end{Definition}

\begin{Remark}
 Taking $\orderedf{k}$ as in the proof of Lemma \ref{sum}, we see
 that our energy can be rewritten as
 \begin{equation*}
  E_{\rho}(f) = \frac{1}{2} \sum_{(x,y)\in \orderedf{1}}
                   \frac{m(x,y)}{|\Gamma_{(x,y)}|} d_Y(f(x),f(y))^2.  
 \end{equation*}

It should be mentioned that 
the above definition of energy is due to Wang \cite{Wan1, Wan2}
when the action of $\Gamma$ is very free.

The energy of a piecewise geodesic map from a finite graph into some
nonpositively curved spaces is considered in \cite{KotSun}, \cite{Leb}.
Our definition of energy generalizes theirs.
\end{Remark}

Note that, in the definition above, we do not assume
that the action of $\Gamma$ is very free.  
If the action of $\Gamma$ is even non-free, one can still construct a 
$\rho$-equivariant map $f:X(0) \longrightarrow Y$ for any $\rho$ as follows:
choose $\mathcal{F}(0)=\{x_1,\dots, x_m\}$ and then $p_i \in Y$
arbitrarily for each $i$. Since the action of $\Gamma$ is properly
discontinuous, the stabilizer $\Gamma_{x_i}$ is finite for each $i$.  
Thus so is $\rho(\Gamma_{x_i})$, and we can find the barycenter $q_i$ of
$\rho(\Gamma_{x_i})p_i$.  Clearly, $q_i$ is fixed by $\rho(\Gamma_{x_i})$. 
Setting $f(x_i)=q_i$ and extending $f$ to whole $X(0)$ by 
$\rho$-equivariance, we obtain a well-defined $\rho$-equivariant map 
$f:X(0)\longrightarrow Y$. 

The assignment $f \mapsto (f(x_1), \dots, f(x_m))$, 
$x_i \in \mathcal{F}(0)$, gives an embedding
of the space of $\rho$-equivariant maps into a product of copies of
$Y$.  Note that $f(x_i)$ must lie in the fixed-point set of
$\rho(\Gamma_{x_i})$, which we denote by $Y_i$.  On the other hand, it is
obvious that any choice of $q_i \in Y_i$ is possible. 
Therefore the space of $\rho$-equivariant maps is identified with the
product space $Y_1 \times \dots \times Y_m$.  Note that
$Y_i$'s are closed convex subsets of $Y$. 

\begin{Definition}
A $\rho$-equivariant map $f : X(0) \longrightarrow Y$ is called
{\em energy-minimizing} if it satisfies $E_\rho(f) \leq E_\rho(g)$
for all $\rho$-equivariant maps $g : X(0) \longrightarrow Y$.
\end{Definition}

We define 
$F_x : (\lk x)(0) \longrightarrow TC_{f(x)}Y$ by 
\begin{equation*}
  F_x(y) = \pi_{f(x)} (f(y)).  
\end{equation*}

\begin{Proposition}
 \label{variation_formula}
  Let $f:X(0) \longrightarrow Y$ be a $\rho$-equivariant map.
If $f$ is energy-minimizing, then
 \begin{equation*}
  \sum_{y \in (\lk x)(0)} m(x,y) \langle W,  F_x(y) \rangle \leq 0
 \end{equation*}
 for all $x \in X(0)$ and $W \in TC_{f(x)}Y$. 
 For each $x \in X(0)$, 
 the barycenter of $\{F_x(y) \mid y \in (\lk x)(0)\}$
 with weight  $\{m(x,y)/m(x) \mid y \in (\lk x)(0)\}$ coincides with
 the origin of $TC_{f(x)}Y$.
\end{Proposition}

\begin{proof} 
Take any $\mathcal{F}(0)$. It suffices to show the assertion of the
proposition for $x\in \mathcal{F}(0)$.
Let $x \in \mathcal{F}(0)$. 
Take a $\rho$-equivariant variation $f_t$ of $f$ so that
$f_t(y)=f(y)$ for $y \in \mathcal{F}(0)$ with $y \not= x$, and that 
$t \mapsto f_t(x)$ is a geodesic starting from $f(x)$. 
 Note that, since $f_t$ is a $\rho$-equivariant variation, $f_t(x)$
 cannot escape from the fixed-point set of $\rho(\Gamma_x)$. 
 Taking $\orderedf{1}=\bigcup_{y \in \mathcal{F}(0)} \orderedf{1}_y$ 
 as in the proof of Lemma \ref{sum}, we see that the part of
 $E_{\rho}(f_t)$ depending on $t$ is 
\begin{equation}
\label{depending_on_t}
  \frac{1}{2}\sum_{(x,y)\in \orderedf{1}_x} 
  \frac{m(x,y)}{|\Gamma_{(x,y)}|} d_Y(f_t(x),f_t(y))^2  
 + \frac{1}{2} \sum_{(y,z) \in \mathcal{F}}
   \frac{m(y,z)}{|\Gamma_{(y,z)}|} 
   d_Y(f_t(y), f_t(z))^2, 
\end{equation}
 where 
 $\mathcal{F}=\{(y,z) \in \orderedf{1} \mid \ z=\gamma x 
                \text{ for some } \gamma \in \Gamma\}$. 
 Consider a map 
 $\Phi: \mathcal{F} \longrightarrow \orderedf{1}_x$ given by
 $(y,z)=(y,\gamma x)\mapsto (x, \gamma^{-1}y)$.  
 Since $\gamma x = \gamma \gamma' x$ if and only if 
 $\gamma' \in \Gamma_x$, $\Phi$ is well-defined, and clearly it is injective.
 It is easy to see that  
 $\Phi(\mathcal{F})=\orderedf{1}_x \setminus \check{X}(1)_x$, where
 $\check{X}(1)_x = \{(x,y) \mid y = \gamma x 
 \text{ for some } \gamma \in \Gamma \}$.
 Hence \eqref{depending_on_t} can be rewritten as
\begin{equation*}
 \begin{split}
& \frac{1}{|\Gamma_x|}\left[
 \sum_{(x,y) \in \ordered{1}_x} m(x,y)
 d_Y(f_t(x),f_t(y))^2 -
 \sum_{(x,y) \in \check{X}(1)_x}\frac{m(x,y)}{2}
 d_Y(f_t(x),f_t(y))^2\right] \\
 = & \frac{1}{|\Gamma_x|}\left[
 \sum_{(x,y) \in \ordered{1}_x\setminus \check{X}(1)_x} m(x,y)
 d_Y(f_t(x),f_t(y))^2 \right. \\
 & \phantom{\frac{1}{|\Gamma_x|}\sum_{(x,y) \in \ordered{1}_x\setminus}}
 \left. + \sum_{(x,y) \in \check{X}(1)_x}\frac{m(x,y)}{2}
 d_Y(f_t(x),f_t(y))^2\right]. 
 \end{split}
\end{equation*}
 Since $f$ is energy-minimizing, we have
\begin{equation}
 \label{minimizing}
 \begin{split}
  &\frac{1}{|\Gamma_x|}\lim_{t \rightarrow +0} \frac{1}{t} 
  \left[\sum_{(x,y) \in \ordered{1}_x\setminus \check{X}(1)_x}m(x,y)
  \left(d_Y(f_t(x),f_t(y))^2 - d_Y(f(x),f(y))^2 \right) \right. \\
 &\left.   + \frac{1}{2}\sum_{(x,y) \in \check{X}(1)_x} m(x,y)
  \left(d_Y(f_t(x),f_t(y))^2-d_Y(f(x),f(y))^2 \right)
  \right] \geq 0 . 
 \end{split}
\end{equation}
 If $(x,y) \in \ordered{1}_x \setminus \check{X}(1)_x$, then
 $f_t(y)=f(y)$ for all $t$.  In this case, we have 
\begin{equation*}
 \begin{split}
   &  \lim_{t \rightarrow +0} 
   \frac{d_Y(f_t(x),f_t(y))^2 - d_Y(f(x),f(y))^2}{t} \\
   = & - 2 d_Y(f(x),f(y)) \cos \angle_{f(x)}(c,f(x)f(y))
 \end{split}
\end{equation*}
 by Lemma \ref{angle}, where $c$ denotes the geodesic $t \mapsto f_t(x)$. 
 On the other hand, if $(x,y) \in \check{X}(1)_x$ and $y=\gamma x$, then 
 we have  
\begin{equation*}
 \begin{split}
 & \lim_{t \rightarrow +0} 
\frac{d_Y(f_t(x), f_t(y))^2- d_Y(f(x),f(y))^2}{t} \\
=& -2d_Y(f(x),f(y))\left[\cos \angle_{f(x)}(f(x)f(y),c)
    + \cos \angle_{f(y)}(f(y)f(x),c')
    \right]
 \end{split}
\end{equation*}
by Lemma \ref{angle} and \ref{angle2}, where $c'$ is a geodesic 
$t \mapsto f_t(y)=\rho(\gamma) f_t(x)$. 
Note that 
$\cos \angle_{f(y)}(f(y)f(x),c')=\cos
 \angle_{f(x)}(f(x)f(\gamma^{-1}x),c)$  
holds, and that  
if $(x,y)=(x,\gamma x) \in \check{X}(1)_x$, then 
$(x,\gamma^{-1}x) \in \check{X}(1)_x$.  Thus we obtain
\begin{equation*}
 \begin{split}
  & \frac{1}{2}\sum_{(x,y) \in \check{X}(1)_x} m(x,y)
     \left(d_Y(f_t(x),f_t(y))^2-d_Y(f(x),f(y))^2\right) \\
  = & -2\sum_{(x,y) \in \check{X}(1)_x} m(x,y) 
     d_Y(f(x),f(y))\cos \angle_{f(x)}(f(x)f(y),c). 
 \end{split}
\end{equation*}
Therefore \eqref{minimizing} becomes
\begin{equation*}
 \begin{split}
 0 &\geq 2 \sum_{y \in (\lk x)(0)} m(x,y) d_Y(f(x),f(y)) \cos
             \angle_{f(x)}(c,f(x)f(y)) \\
   &= 2 \sum_{y \in (\lk x)(0)} m(x,y) \langle V, F_x(y) \rangle, 
 \end{split}
\end{equation*}
 where $V=([c],1) \in TC_{f(x)}Y$ and $[c]$ is the equivalence class of
 $c$ in $S_{f(x)}Y$. By the continuity of the inner product and Lemma
 \ref{fixed-point_set}, we see that  
\begin{equation}
\label{barycenter_inequality}
    0 \geq  \sum_{y \in (\lk x)(0)} m(x,y) \langle V, F_x(y) \rangle
\end{equation}
 holds for all $V \in (TC_{f(x)}Y)^{\rho(\Gamma_x)}$. 
 Let $V=(V_0,t) \in (TC_{f(x)}Y)^{\rho(\Gamma_x)}$, where $V_0=V/|V|$
 and $t=|V|$. 
 By the definition of the distance on the tangent
 cone and \eqref{barycenter_inequality}, we see that  
\begin{equation}\label{barycenter_on_cone}
 \begin{split}
 & \sum_{y \in (\lk x)(0)} m(x,y)
   d_{TC_{f(x)}Y}((V_0,t),F_x(y))^2 \\
 =& \sum_{y \in (\lk x)(0)} m(x,y)\big[t^2 + |F_x(y)|^2  
    - 2t\langle V_0,F_x(y)\rangle \big]
 \end{split}
\end{equation}
 is an increasing function of $t \geq 0$.  
 In other words,  the function
 on  $(TC_{f(x)}Y)^{\rho(\Gamma_x)}$ with variable $V=(V_0,t)$ defined
 by the left-hand side of \eqref{barycenter_on_cone} is minimized at the 
 origin of $TC_{f(x)}Y$. Recall that the barycenter of 
 $\{F_x(y)\mid y \in (\lk x)(0)\}$ with weight 
 $\{m(x,y)/m(x) \mid y \in (\lk x)(0)\}$ is the point that minimizes
 the left-hand side of \eqref{barycenter_on_cone} 
 (divided by the constant $m(x)$) on the whole $TC_{f(x)}Y$.  
 Since $\{F_x(y)\mid y \in (\lk x)(0)\}$ (resp.~$m(x,y)$) is 
 $\rho(\Gamma_x)$-invariant (resp.~$\Gamma$-invariant), their barycenter  
 must lie in $(TC_{f(x)}Y)^{\rho(\Gamma_x)}$.  
 Therefore the origin of $TC_{f(x)}Y$ must be the barycenter
 of $\{F_x(y)\mid y \in (\lk x)(0)\}$.  
 Noting that the origin of $TC_{f(x)}Y$ minimizes the left-hand side of
 \eqref{barycenter_on_cone}, we see that the inequality in our
 proposition follows from computation similar to (but substantially easier
 than) the one we have done in order to deduce 
 \eqref{barycenter_inequality}. 
\end{proof}

\begin{Definition}
Let $- \Delta f(x) \in TC_{f(x)}Y$ be the barycenter of 
$\{F_x(y)\mid y \in (\lk x)(0)\}$ 
with weight $\{m(x,y)/m(x) \mid y \in (\lk x)(0)\}$.  
We say that a $\rho$-equivariant map $f:X(0)\longrightarrow Y$
is {\it harmonic} if $- \Delta f(x) = 0_{f(x)}$ for all 
$x \in \mathcal{F}(0)$, where $0_{f(x)}$ is the origin 
of $TC_{f(x)}Y$.
\end{Definition}

\begin{Remark}
Gromov \cite{Gro2} defines the harmonicity of an equivariant map
in a slightly different way, without referring to the tangent cone.
It turns out, however, that the two definitions are equivalent.
\end{Remark} 

When the tangent cone $TC_{f(x)}Y$ is isometric to a Hilbert space, 
it is obvious that 
\begin{equation*}
- \Delta f(x)= \sum_{y \in (\lk x)(0)} \frac{m(x,y)}{m(x)}F_x(y)
\end{equation*}
holds.  
In general, the right-hand side does not make sense, as the addition is 
not defined on a non-Hilbertian tangent cone.
However, we still have the following formula.
\begin{Lemma}
\label{Delta_f}
 Let $\xi_0(x)=0_{f(x)}$ if $- \Delta f(x)=0_{f(x)}$, and  
 $\xi_0(x)= - \Delta f(x)/ |- \Delta f(x)| \in S_{f(x)}Y$ otherwise. 
 Then we have 
 \begin{equation*}
    |- \Delta f(x)|= \sum_{y \in (\lk x)(0)} 
     \left\langle \xi_0(x), \frac{m(x,y)}{m(x)}F_x(y) \right\rangle. 
 \end{equation*} 
\end{Lemma}

\begin{proof}
 If $- \Delta f(x)=0_{f(x)}$, the formula is obvious. 
 Suppose $- \Delta f(x) \not= 0_{f(x)}$, and set $\eta=t\xi_0(x)$. 
 Consider a function of $t$ defined by
 \begin{equation}
  \label{barycenter_functional}
  \begin{split}
       &  \sum_{y \in (\lk x)(0)}\frac{m(x,y)}{m(x)} 
          d_{C_{f(x)}Y}(F_x(y), \eta)^2 \\
       = &  \sum_{y \in (\lk x)(0)} \frac{m(x,y)}{m(x)} \left(
         |F_x(y)|^2 + t^2 - 
         2t\langle \xi_0(x), F_x(y)\rangle \right). 
  \end{split}
 \end{equation}
 The right-hand side of \eqref{barycenter_functional} is
 minimized at 
 \begin{equation*}
  t = \sum_{y \in (\lk x)(0)} \frac{m(x,y)}{m(x)}
  \langle \xi_0(x), F_x(y)\rangle. 
 \end{equation*}
 Since $- \Delta f(x)= |- \Delta f(x)|\xi_0$ is the point that minimizes the
 left-hand side of  \eqref{barycenter_functional} regarded as a function
 defined for all $\eta \in TC_{f(x)}Y$, this $t$ must be equal to 
 $|- \Delta f(x)|$.   
\end{proof}

Wang \cite{Wan2} proved an existence result for an energy-minimizing 
$\rho$-equivariant map.
His argument extends to our setting without difficulty and gives
the following 

\begin{Theorem}\label{existence}
Let $X$ be a simplicial complex equipped
with an admissible weight, and let $\Gamma$ be a finitely generated
group acting by automorphisms, properly discontinuously and cofinitely
on $X$.
Let $Y$ be a Hadamard space, and $\rho : \Gamma \longrightarrow \isom(Y)$
a homomorphism.
Suppose that $Y$ is locally compact and $\rho$ is reductive 
in the sense of Jost $($see \cite{Jos3}$)$.
Then there exists an energy-minimizing $\rho$-equivariant map
$f : X(0) \longrightarrow Y$.
\end{Theorem}

\section{Bochner-type formulas for equivariant maps}

Let $X$ and $\Gamma$ be as in the preceding section.
Let $Y$ be a Hadamard space and $\rho : \Gamma
\longrightarrow \isom(Y)$ a homomorphism.
We shall write down Bochner-type formulas for a $\rho$-equivariant
map $f : X(0) \longrightarrow Y$.

\begin{Proposition}\label{nonlinear_case}
Let $f : X(0) \longrightarrow Y$ be a $\rho$-equivariant map.
Then the following two formulas hold{\rm :}
\begin{eqnarray}\label{wformula1}
0 &=& \sum_{x\in \mathcal{F}(0)} \frac{1}{|\Gamma_x|}
\left[ \sum_{(y,y') \in (\overrightarrow{\lk x})(1)}
m(x,y,y')\, d_{TC_{f(x)}Y}(F_x(y), F_x(y'))^2 \right. \nonumber\\
&& - \sum_{y\in (\lk x)(0)} m(x,y)\, d_{TC_{f(x)}Y}(0_{f(x)},
F_x(y))^2 \\
&& \left. + \sum_{(y,y') \in (\overrightarrow{\lk x})(1)}
m(x,y,y') \left\{ d_Y(f(y), f(y'))^2 - d_{TC_{f(x)}Y}(F_x(y),
F_x(y'))^2 \right\} \right];\nonumber
\end{eqnarray}
\begin{eqnarray}\label{wformula1'}
\lefteqn{\sum_{x\in \mathcal{F}(0)} \frac{m(x)}{|\Gamma_x|}
d_{TC_{f(x)}Y}(0_{f(x)},- \Delta f(x))^2} \nonumber \\
&=& \sum_{x\in \mathcal{F}(0)} \frac{1}{|\Gamma_x|}
\left[ \sum_{(y,y') \in (\overrightarrow{\lk x})(1)}
m(x,y,y')\, d_{TC_{f(x)}Y}(F_x(y), F_x(y'))^2 \right. \nonumber\\
&& - \sum_{y\in (\lk x)(0)} m(x,y)\, d_{TC_{f(x)}Y}(- \Delta f(x),
F_x(y))^2 \\
&& \left. + \sum_{(y,y') \in (\overrightarrow{\lk x})(1)}
m(x,y,y') \left\{ d_Y(f(y), f(y'))^2 - d_{TC_{f(x)}Y}(F_x(y),
F_x(y'))^2 \right\} \right].\nonumber
\end{eqnarray}
\end{Proposition}

\begin{proof}
We shall regard the correspondence $\overrightarrow{X}(1) \ni (x,y)
\mapsto F_x(y) \in TC_{f(x)}Y$ as an analogue of equivariant
$1$-cocycle (see Appendix). 
Note that \eqref{dL2} (with $\alpha$ a $1$-cocycle) 
can be rewritten as
\begin{eqnarray*}
0 &=& ||d\alpha||_{L^2}^2 \\
&=& \sum_{x\in \mathcal{F}(0)} \frac{1}{|\Gamma_x|}
\left[ \frac{1}{2}\sum_{y\in (\lk x)(0)} m(x,y) ||\alpha(x,y)||^2
\right. \\
&& \left. - \sum_{(y,y')\in (\overrightarrow{\lk x})(1)}
m(x,y,y') \langle \alpha(x,y), \alpha(x,y') \rangle \right] \\
&=& \frac{1}{2} \sum_{x\in \mathcal{F}(0)} \frac{1}{|\Gamma_x|}
\left[ \sum_{(y,y') \in (\overrightarrow{\lk x})(1)}
m(x,y,y') ||\alpha(x,y)-\alpha(x,y')||^2 \right. \\
&& \left. - \sum_{y\in (\lk x)(0)} m(x,y)
||\alpha(x,y)||^2 \right].
\end{eqnarray*}
The analogue of the last expression in the present setting is
\begin{eqnarray}\label{eq_a}
&&\frac{1}{2} \sum_{x\in \mathcal{F}(0)} \frac{1}{|\Gamma_x|}
\left[ \sum_{(y,y') \in (\overrightarrow{\lk x})(1)}
m(x,y,y')\, d_{TC_{f(x)}Y}(F_x(y), F_x(y'))^2 \right. \\
&& \left. - \sum_{y\in (\lk x)(0)} m(x,y)\,
d_{TC_{f(x)}Y}(0_{f(x)}, F_x(y))^2 \right]. \nonumber
\end{eqnarray}
This expression is not equal to zero in general; there are extra terms
reflecting the curvature of $Y$, and we shall identify them.
First note that
$d_{TC_{f(x)}Y}(0_{f(x)}, F_x(y)) = d_Y(f(x), f(y))$.
Using the defining properties of the admissible weight $m$ and
Lemma \ref{sum}, we compute
\begin{eqnarray*}
&& \sum_{x\in \mathcal{F}(0)} \frac{1}{|\Gamma_x|}
\sum_{y\in (\lk x)(0)} m(x,y)\, d_{TC_{f(x)}Y}(0_{f(x)},
F_x(y))^2 \\
&=& \sum_{x\in \mathcal{F}(0)} \frac{1}{|\Gamma_x|}
\sum_{(x,y,y')\in \overrightarrow{X}(2)_x} m(x,y,y')\,
d_Y(f(x), f(y))^2\\
&=& \sum_{(x,y,y')\in \overrightarrow{\mathcal{F}}(2)}
\frac{m(x,y,y')}{|\Gamma_{(x,y,y')}|}\, d_Y(f(x), f(y))^2 \\
&=& \sum_{(x,y,y')\in \overrightarrow{\mathcal{F}}(2)}
\frac{m(x,y,y')}{|\Gamma_{(x,y,y')}|}\, d_Y(f(y), f(y'))^2 \\
&=& \sum_{x\in \mathcal{F}(0)} \frac{1}{|\Gamma_x|}
\sum_{(y,y')\in (\overrightarrow{\lk x})(1)} m(x,y,y')\,
d_Y(f(y), f(y'))^2.
\end{eqnarray*}
For the third equality, we have used the argument as in the proof
of Proposition \ref{adjoint}.
It follows that
\begin{eqnarray*}
(\ref{eq_a}) &=&\frac{1}{2} \sum_{x\in \mathcal{F}(0)}
\frac{1}{|\Gamma_x|} \sum_{(y,y') \in (\overrightarrow{\lk x})(1)}
m(x,y,y')\, \bigl\{ d_{TC_{f(x)}Y}(F_x(y), F_x(y'))^2 \\
&& - d_Y(f(y), f(y'))^2 \bigr\}.
\end{eqnarray*}
Rearranging this gives (\ref{wformula1}).

To prove \eqref{wformula1'}, we compute
\begin{eqnarray*}
\lefteqn{\sum_{y\in (\lk x)(0)} m(x,y)\, d_{TC_{f(x)}Y}(- \Delta f(x),
F_x(y))^2} \\
&=& \sum_{y\in (\lk x)(0)} m(x,y) \left( |- \Delta f(x)|^2
+ |F_x(y)|^2 - 2 \langle - \Delta f(x), F_x(y) \rangle \right) \\
&=& m(x)\, |- \Delta f(x)|^2 + \sum_{y\in (\lk x)(0)} m(x,y)\,
|F_x(y)|^2 - 2 m(x)\, |- \Delta f(x)|^2 \\
&=& - m(x)\, d_{TC_{f(x)}Y}(0_{f(x)}, - \Delta f(x))^2
+ \sum_{y\in (\lk x)(0)} m(x,y)\, d_{TC_{f(x)}Y}(0_{f(x)}, F_x(y))^2.
\end{eqnarray*}
For the second equality, we have used Lemma \ref{Delta_f}. 
\eqref{wformula1'} now follows from \eqref{wformula1}.
The proof of Proposition \ref{nonlinear_case} is completed.
\end{proof}

\begin{Remark}\label{3rd_term}
Since the map $\pi_{f(x)} : Y \longrightarrow TC_{f(x)}Y$
is distance non-increasing, we have
$$
d_Y(f(y), f(y')) \geq d_{TC_{f(x)}Y}(F_x(y), F_x(y')),\quad
(y, y') \in (\overrightarrow{\lk x})(1).
$$
Hence the third terms in $[\quad]$ of \eqref{wformula1},
\eqref{wformula1'} are nonnegative.
\end{Remark}

\vskip0.5cm
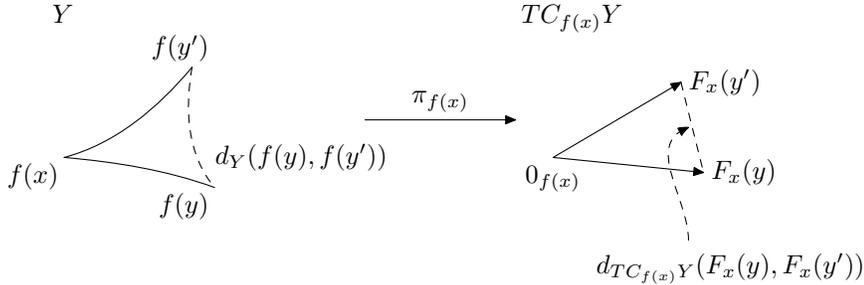
\hspace{.5cm}
\includegraphics{nonpositive.1}
\vspace{.5cm}

We recall the definition of the numerical invariant
introduced by Wang \cite{Wan2}.
\begin{Definition}
For $x\in X(0)$ and a Hadamard space $T$ (it is typically 
a tangent cone of a Hadamard space), define
$$
\lambda_1(\lk x, T) = \inf_\varphi \frac{
\frac{1}{2}\sum_{(y,z)\in (\overrightarrow{\lk x})(1)} m(x,y,z)\,
d_T(\varphi(y), \varphi(z))^2}{
\sum_{y\in (\lk x)(0)} m(x,y)\, d_T(\overline{\varphi}, 
\varphi(y))^2},
$$
where the infimum is taken over all nonconstant maps 
$\varphi : (\lk x)(0) \longrightarrow T$ and $\overline{\varphi}$ 
is the barycenter of $\{\varphi(y) \mid y\in (\lk x)(0)\}$ 
with weight $\{m(x,y)/m(x) \mid y\in (\lk x)(0)\}$.
\end{Definition}

Wang \cite[Theorem 7.3]{Wan2} proved the following result, assuming 
$\Gamma\backslash X$ was a simplicial complex.
Suppose that $Y$ is locally compact and $\rho$ is reductive.
If $\lambda_1(\lk x, TC_pY) > 1/2$ for all $x\in X(0)$ and $p\in Y$, 
then $\rho(\Gamma)$ has a fixed point in $Y$.
Indeed, by Theorem \ref{existence}, there exists an energy-minimizing
$\rho$-equivariant map $f : X(0) \longrightarrow Y$.
Then it follows from either of the formulas of Proposition 
\ref{nonlinear_case} and Remark \ref{3rd_term} that $f$ must be 
a constant map, and thus the image point of $f$ is fixed by 
$\rho(\Gamma)$.

\section{Gradient flow of the energy functional and fixed-point
theorems} 

In this section, we shall strengthen Wang's fixed-point theorem mentioned
in the preceding section by deriving the same conclusion without assuming
$Y$ is locally compact nor $\rho$ is reductive.
Our proof uses the gradient flow, introduced by Jost \cite{Jos4} and Mayer \cite{May}, 
associated with our energy functional.

Let $X$ be a simplicial complex equipped with an admissible weight $m$. 
We assume $\lk x$ is connected for all $x \in X(0)$. 
Let $\Gamma$ be a finitely generated group acting by
automorphisms, properly discontinuously and cofinitely on $X$. Let $Y$
be a Hadamard space, and let
$\rho:\Gamma \longrightarrow \isom (Y)$ be a homomorphism. 
As we have seen in \S 2, the space of $\rho$-equivariant maps
$\mathcal{M}$ can be identified with a product space 
$\prod_{x \in \mathcal{F}(0)}Y_x$, where $Y_x \subset Y$ is the
fixed-point set of $\rho(\Gamma_x)$.  We define a metric on $\mathcal{M}$
by 
\begin{equation*}
d_{\mathcal{M}}(f_0,f_1)^2 =\sum_{x \in \mathcal{F}(0)}
\frac{m(x)}{|\Gamma_x|}d_Y(f_0(x),f_1(x))^2, \quad f_0,f_1 
\in \mathcal{M}. 
\end{equation*}
Since $Y_x$'s are closed convex subsets of the Hadamard space $Y$,
$(\mathcal{M},d_{\mathcal{M}})$ is again a Hadamard space. 
It is clear
that $E_{\rho}$ is continuous on $(\mathcal{M},d_{\mathcal{M}})$. 
The inner product of $TC_{f}\mathcal{M}$ takes the form of 
$\sum_{x \in \mathcal{F}(0)}(m(x)/|\Gamma_x|)
\langle \cdot , \cdot \rangle_{TC_{f(x)}Y_x}$. 

Let $\mathcal{F}(0)=\{x_1, \dots, x_m\}$. 
Then, with respect to $d_{\mathcal{M}}$, the geodesic $c$ in $\mathcal{M}$
joining $f_0$ and $f_1$ is given by  
\begin{equation*}
c(t) = (c_{x}(d_x t/ d))_{x \in \mathcal{F}(0)}
=(c_{x_1}(d_{x_1} t/ d), \dots , c_{x_m}(d_{x_m} t/d)),
\end{equation*}
where $c_{x}$ is the geodesic in $Y_x$ joining $f_0(x)$ and $f_1(x)$, 
$d_x = d_Y(f_0(x), f_1(x))$, and $d=d_{\mathcal{M}}(f_0,f_1)$. 
By \cite[Proposition 2.2]{BriHae}, it is readily seen that
\begin{equation*}
 d_Y(c(tT),c'(tT'))^2 \leq (1-t)d_Y(c(0),c'(0))^2 + 
                           td_Y(c(T),c'(T'))^2
\end{equation*}
holds for any pair of geodesics 
$c:[0,T]\longrightarrow Y$ and $c':[0,T']\longrightarrow Y$.
Therefore $E_{\rho}$ is a convex function on $\mathcal{M}$, that is,
$E_{\rho}$ satisfies
\begin{equation*}
 E_{\rho}(f_t)\leq (1-t)E_{\rho}(f_0) + t E_{\rho}(f_1), 
\end{equation*}
where $f_t \in \mathcal{M}$ is the point with a fraction $t$ from $f_0$
to $f_1$. 

A slight modification of the proof of Proposition
\ref{variation_formula} shows the following proposition. 

\begin{Proposition}\label{variation_formula2}
 Let $c(t)=(c_x(d_xt/d))_{x \in \mathcal{F}(0)}$ be a geodesic in
 $\mathcal{M}$ starting from $f$ expressed as above, 
 and $[c]$ the element of the 
 space of directions at $f=c(0) \in \mathcal{M}$ defined by $c$, namely,  
 $[c]=(d_x[c_x]/d)_{x \in \mathcal{F}(0)} \in TC_{f}\mathcal{M}$. 
 Then we have  
  \begin{eqnarray}\label{variation_f}
\lefteqn{\lim_{t \rightarrow +0}\frac{E_{\rho}(c(t))-E_{\rho}(c(0))}{t}}\\
&=& -\sum_{x \in \mathcal{F}(0)}\frac{m(x)}{|\Gamma_x|} 
\sum_{y \in (\lk x)(0)}
\left\langle \frac{d_x}{d}[c_x], 2\, \frac{m(x,y)}{m(x)}F_x(y)
\right\rangle_{TC_{f(x)}Y}. \nonumber
\end{eqnarray}
\end{Proposition}
Note that if $Y$ is a Riemannian manifold, the right-hand side of 
\eqref{variation_f} can be rewritten as 
$-\left\langle [c], (2(- \Delta f(x)))_{x \in \mathcal{F}(0)}
\right\rangle_{\mathcal{M}}$.

Let $\xi_0(x)=- \Delta f(x)/|- \Delta f(x)|$ as in \S 2. 
Suppose that for every $x \in \mathcal{F}(0)$ there exists a geodesic $c_x$ 
satisfying $\xi_0(x)= [c_x]$, and consider the geodesic $c$ in $\mathcal{M}$ 
defined by 
\begin{equation*}
c(t)=\left(c_x \left( \alpha_x t \right)\right)_{x \in \mathcal{F}(0)},
\quad \alpha_x = \frac{|- \Delta f(x)|}{|(- \Delta f(x))_{x \in \mathcal{F}(0)}|}.
\end{equation*}
Then, by the proposition above and Lemma \ref{Delta_f}, we have 
\begin{equation*}
 \lim_{t \rightarrow +0} \frac{E_{\rho}(c(t))-E_{\rho}(c(0))}{t}
 =|(2(- \Delta f(x)))_{x \in \mathcal{F}(0)}|.
\end{equation*}
Though such a geodesic $c_x$ may not exist in general,
we still have a sequence of geodesics $\{c_j\}_{j=1}^{\infty}$ in $\mathcal{M}$
satisfying
\begin{equation}
\label{gradient_approximated}
 \lim_{j \rightarrow \infty} 
 \lim_{t \rightarrow +0} \frac{E_{\rho}(c_j(t))-E_{\rho}(c_j(0))}{t}
 =|(2(- \Delta f(x)))_{x \in \mathcal{F}(0)}|.
\end{equation}
Motivated by Proposition \ref{variation_formula2} and the observation 
we have just made, we set
\begin{equation*}
 -\mathrm{grad}E_{\rho}(f)
    = (2(- \Delta f(x)))_{x \in \mathcal{F}(0)}. 
\end{equation*}

Jost \cite{Jos4} and Mayer \cite{May} defined a gradient flow for a certain class 
of functions 
on a nonpositively curved metric space.  Their theory is applicable to our
convex function $E_{\rho}$ on the Hadamard space $\mathcal{M}$.  
We denote the gradient flow of $E_{\rho}$ starting from $f$ by $f(t)$. 
Then $f(t)$ is defined for all $t \in [0,\infty)$ 
(\cite[Theorem 1.13]{May}).  As one expects, 
$t \mapsto E_{\rho}(f(t))$ is a non-increasing continuous
function, and moreover, 
$E_{\rho}(f(t))$ is Lipschitz continuous on each closed interval
$[t,t']\subset (0,\infty)$ (\cite[Corollary 2.11]{May}). 
Though Mayer did not define a gradient vector,
he introduced ``the norm of the gradient vector''
$|\nabla_{-}E_{\rho}|$ defined by 
\begin{equation*}
 |\nabla_{-}E_{\rho}|(f_0)=\max
 \left\{ \limsup_{f \rightarrow f_0}
 \frac{E_{\rho}(f_0)-E_{\rho}(f)}{d_{\mathcal{M}}(f_0,f)}, 
 0  \right\}. 
\end{equation*}
Mayer showed that this function has some properties satisfied by the norm of
the usual gradient vector, such as 
\begin{eqnarray}
  & \displaystyle{
     \lim_{s \rightarrow +0} \frac{d_{\mathcal{M}}(f(t+s),f(t))}{s}
  =|\nabla_{-}E_{\rho}|(f(t)) \quad \text{for all $t$} }, \notag \\
  & \label{property2}
    \displaystyle{
     - \frac{dE_{\rho}(f(t))}{dt}=
   \left( |\nabla_{-}E_{\rho}|(f(t))\right)^2 \quad
   \text{for almost all $t>0$}. } 
\end{eqnarray}
See \cite[Theorem 2.17, Corollary 2.18]{May}.
Since $E_{\rho}$ is a convex function, $|\nabla_{-}E_{\rho}|$ satisfies 
\begin{equation}
\label{property3}
    \displaystyle{  
    |\nabla_{-}E_{\rho}|(f(t))=\sup_{s>0} 
    \frac{d_{\mathcal{M}}(f(t),f(t+s))}{s} \quad \text{for all }t }
\end{equation}
(\cite[Proposition 2.34]{May}). He also proved that 
$t \mapsto |\nabla_{-}E_{\rho}|f(t)$ is right continuous
(\cite[Corollary 2.28]{May}).   Note that we have 
\begin{equation}
\label{grad_vs_nabla}
 |\nabla_{-}E_{\rho}|(f) \geq |-\mathrm{grad}E_{\rho}|(f)
\end{equation}
by 
\eqref{gradient_approximated} and the definition of 
$|\nabla_{-}E_{\rho}|$.   

\begin{Remark}
It is plausible that the inequality \eqref{grad_vs_nabla} is 
actually an
equality, and that $-\mathrm{grad}E_{\rho}(f(t))$ gives the velocity vector 
of the gradient flow $f(t)$, though we will not prove these statements
in this paper. 
\end{Remark}

Our Bochner-type formula \eqref{wformula1'} gives the following
estimate of $|\nabla_{-} E_{\rho}|$.

\begin{Lemma}
\label{gradient_estimate}
Suppose that there exists a constant $C > 1/2$ such that 
$\lambda_1(\lk x, TC_pY)\geq C$ for all $x \in X(0)$ and $p \in Y$.  
Then 
\begin{equation*}
|\nabla_{-}E_{\rho}|^2(f) \geq 4 C E_{\rho}(f)
\end{equation*}
holds for all $f \in \mathcal{M}$. 
\end{Lemma}

\begin{Remark}
By employing \.Zuk's idea \cite{Zuk1}, the assumption on Wang's invariant
in the lemma, hence in Theorem \ref{fixed_point_theorem1} below, can be 
weakened to the following:
there exists a constant $C' > 1$ such that
$$
\lambda_1(\lk x, TC_pY) + \lambda_1(\lk y, TC_qY) \geq C' \quad
\mbox{for all $\{x, y\} \in X(1)$ and $p,q \in Y$}.
$$
\end{Remark}

\begin{proof}
We prove the assertion of the lemma under the assumption of the remark.
Let $f\in \mathcal{M}$, and write $\lambda(x) = \lambda_1(\lk x, TC_{f(x)}Y)$
for simplicity.
By \eqref{wformula1'} and Remark \ref{3rd_term}, it is immediate that
\begin{eqnarray*}
\lefteqn{\frac{1}{4} |-\mathrm{grad }E_{\rho}(f)|^2}\\
&\geq& \sum_{x\in \mathcal{F}(0)} \frac{1}{|\Gamma_x|}
(2 \lambda(x) - 1)
\sum_{y\in (\lk x)(0)} m(x,y)\, d_{TC_{f(x)}Y}(- \Delta f(x),
F_x(y))^2.
\end{eqnarray*}
By using the computation in the proof of \eqref{wformula1'}, we
estimate the right-hand side as follows:
\begin{eqnarray*}
{\rm R.H.S.} &=& \sum_{x\in \mathcal{F}(0)} \frac{1}{|\Gamma_x|}
(2 \lambda(x) - 1)
\biggl( - m(x)\,d_{TC_{f(x)}Y}(0_{f(x)}, - \Delta f(x))^2\\
&& 
+ \sum_{y\in (\lk x)(0)} m(x,y)\,d_{TC_{f(x)}Y}(0_{f(x)}, F_x(y))^2
\biggr) \\
&\geq& - 3 \sum_{x\in \mathcal{F}(0)} \frac{m(x)}{|\Gamma_x|}
d_{TC_{f(x)}Y}(0_{f(x)}, - \Delta f(x))^2\\
&& + \sum_{(x,y)\in \orderedf{1}} \frac{m(x,y)}{|\Gamma_{(x,y)}|}
(2 \lambda(x)- 1) d_Y(f(x), f(y))^2,
\end{eqnarray*}
where we have used the fact that $\lambda(x) \leq 2$; see Proposition
\ref{wang_invariant} below.
Therefore, we obtain (c.f.~the proof of Proposition 
\ref{zuk's_formulation})
\begin{eqnarray*}
|-\mathrm{grad }E_{\rho}(f)|^2 &\geq& 
\sum_{(x,y)\in \orderedf{1}} \frac{m(x,y)}{|\Gamma_{(x,y)}|}
(\lambda(x) + \lambda(y) - 1) 
d_Y(f(x), f(y))^2\\
&\geq& 2 C' E_{\rho}(f).
\end{eqnarray*}
Now \eqref{grad_vs_nabla} gives the desired inequality. 
\end{proof}

Lemma \ref{gradient_estimate} enables us to prove the following
fixed-point theorem, which is a strengthened version of 
Wang's result mentioned in \S 3.

\begin{Theorem}
\label{fixed_point_theorem1}
 Let $X$ be a simplicial complex equipped with an
 admissible weight $m$. 
 We assume $\lk x$ is connected for all  $x \in X(0)$. 
 Let $\Gamma$ be a finitely generated group acting by
 automorphisms, properly discontinuously and cofinitely on $X$. Let $Y$
 be a Hadamard space.
Suppose that there exists a constant $C > 1/2$ such that
$\lambda_1(\lk x, TC_pY)\geq C$ for all $x \in X(0)$ and $p \in Y$.
Then for any homomorphism $\rho:\Gamma \longrightarrow \isom (Y)$,
$\rho(\Gamma)$ has a fixed point in $Y$. 
\end{Theorem}

\begin{proof}
Take any $f \in \mathcal{M}$.  We shall show that the gradient flow
$f(t)$ starting from $f$ converges to a constant map as 
$t \rightarrow \infty$.  Let $E(t)=E_{\rho}(f(t))$ and 
$|\nabla_{-}E|(t)=|\nabla_{-}E_{\rho}|(f(t))$. Set $E(t)=e^{-h(t)}$. 
Since $E(t)$ is a locally Lipschitz continuous function on $(0,\infty)$, 
$E'(t)=-h'(t)e^{-h(t)}$ exists for almost all $t$.
Hence, by \eqref{property2} and Lemma \ref{gradient_estimate}, 
$h'(t) \geq C'$ for almost all $t$. 
Since $h(t)$ is locally Lipschitz continuous, there is a constant $C''$ 
such that $h(t) \geq C't + C''$ for all $t$. 

 Take any closed interval $[t,t']\subset (0,\infty)$ and set 
 \begin{equation*}
 s_j^{(n)}=t+\frac{t'-t}{2^n} j, \quad j=0, \dots, 2^n.
 \end{equation*}
 Then, by \eqref{property3}, we obtain
 \begin{equation*}
 d_{\mathcal{M}}(f(t),f(t')) \leq
 \sum_{j=1}^{2^n}d_{\mathcal{M}}(f(s_{j-1}^{(n)}),f(s_j^{(n)})) \leq
 \sum_{j=1}^{2^n} \frac{t'-t}{2^n} |\nabla_{-}E|(s_{j-1}^{(n)}). 
\end{equation*} 
 Define a simple function $F:[t,t']\longrightarrow \R$ by
 \begin{equation*}
 F_n(s)= |\nabla_{-}E|(s_j^{(n)}), 
 \quad s \in (s_{j-1}^{(n)},s_j^{(n)}],
 \end{equation*}
so that
 \begin{equation*}
 \int_{t}^{t'}F_n(s)\ ds =
 \sum_{j=1}^{2^n} \frac{t'-t}{2^n} |\nabla_{-}E|(s_{j}^{(n)}).
 \end{equation*}
Since $E(s)$ is Lipschitz continuous on $[t,t']$ and \eqref{property2}
holds for almost all $s$, we see that $|\nabla_{-}E|(s)$ is 
essentially bounded on $[t,t']$.
Then, since $|\nabla_{-}E|(s)$ exists for all $s$ and is right continuous, 
$|\nabla_{-} E|(s)$ is bounded on $[t,t']$. 
Therefore, taking a subsequence if necessary, the above integral
converges as $n\rightarrow \infty$.
In particular, we obtain
 \begin{equation*}
  d_{\mathcal{M}}(f(t),f(t')) \leq 
  \lim_{n \rightarrow \infty} \int_{t}^{t'}F_n(s) ds. 
 \end{equation*}
 On the other hand, the right continuity of $|\nabla_{-}E|(s)$
 guarantees the pointwise convergence of $\{F_n\}_{n=1}^{\infty}$ to 
 $|\nabla_{-}E|(s)$ as $n \rightarrow \infty$.  Therefore, by the
 dominated convergence theorem, 
 \begin{equation*}
 \lim_{n \rightarrow \infty} \int_{t}^{t'}F_n(s) ds
  = \int_{t}^{t'}|\nabla_{-} E|(s) ds. 
 \end{equation*}

Since $|\nabla_{-}E|(s)=\sqrt{-E'(s)}$ for almost all $s$, we see that
\begin{eqnarray*}
d_{\mathcal{M}}(f(t),f(t')) &\leq& \int_{t}^{t'}\sqrt{-E'(s)} ds 
= \int_{t}^{t'}\sqrt{h'(s)}e^{-h(s)/2} ds \\
&\leq& \int_{t}^{t'}\frac{h'(s)}{\sqrt{C'}}e^{-h(s)/2} ds 
\leq \frac{2}{\sqrt{C'}} e^{-h(t)/2} \leq C''' e^{-C't/2}.
\end{eqnarray*}
Note that the second inequality follows from $h'(t) \geq C'$. 
Therefore, for any divergent sequence $\{t_j\}_{j=1}^{\infty} \subset (0,\infty)$,
 $\{f(t_j)\}_{j=1}^{\infty}$ is a Cauchy sequence in $\mathcal{M}$, and
it has a limit $f_{\infty}$.
 Since $E_{\rho}(f(t)) \leq e^{-Ct-C'}$, $E_{\rho}(f_{\infty})=0$,
that is, $f_{\infty}$ is a constant map.  This completes the proof. 
\end{proof}

\section{Computation of $\lambda_1$}

In this section, we shall give lower and upper estimates of the 
invariant $\lambda_1$.
We shall also compute $\lambda_1$ in some cases, and deduce more 
concrete fixed-point theorems from Theorem \ref{fixed_point_theorem1}.

Let $X$ be a simplicial complex equipped with an admissible weight $m$.
Throughout this section, we shall assume that $\lk x$ is connected 
for all $x\in X(0)$.
Let $\mathcal{H}$ be a Hilbert space with inner product 
$\langle \cdot, \cdot \rangle$, and for $l=0,1$,
let $C^l(\lk x, \mathcal{H})$ denote the set of simplicial $l$-cochains 
on $\lk x$ with values in $\mathcal{H}$.
We define inner products and norms on $C^l(\lk x, \mathcal{H})$,
$l=0,1$, by
$$
\langle f,g \rangle_{L^2} = \sum_{y\in (\lk x)(0)}m(x,y)
\langle f(y), g(y) \rangle,\quad f,g\in C^0(\lk x, \mathcal{H});
$$
$$
\langle \alpha, \beta \rangle_{L^2} = \frac{1}{2} \sum_{(y,y') \in
(\overrightarrow{\lk x})(1)}m(x,y,y')
\langle \alpha(y,y'), \beta(y, y') \rangle,\quad \alpha,\beta \in 
C^1(\lk x, \mathcal{H});
$$
$$
||f||_{L^2}^2 = \langle f,f \rangle_{L^2};\quad 
||\alpha||_{L^2}^2 = \langle \alpha, \alpha \rangle_{L^2}.
$$
We set $\mu_1(\lk x) = \lambda_1(\lk x, \mathcal{H})$, namely,
$$
\mu_1(\lk x) = \inf \frac{||df||_{L^2}^2}{||f-\overline{f}||_{L^2}^2},
$$
where the infimum is taken over all nonconstant 
$f\in C^0(\lk x, \mathcal{H})$ and 
$$
\overline{f} = \sum_{y\in (\lk x)(0)} \frac{m(x,y)}{m(x)} f(y).
$$
This quantity coincides with the first nonzero eigenvalue of 
the Laplacian
$\Delta_{\lk x}: C^0(\lk x, \mathcal{H}) \longrightarrow C^0(\lk x, 
\mathcal{H})$, given by
$$
(\Delta_{\lk x} f) (y) = f(y) - 
\sum_{y' ; \, (y,y')\in (\overrightarrow{\lk x})(1)} 
\frac{m(x,y,y')}{m(x,y)} f(y').
$$
It is an easy matter to verify that $\mu_1(\lk x)$ does not depend 
on the choice of $\mathcal{H}$.

Let $T$ be a Hadamard space.
In Proposition \ref{wang_invariant} below, we shall 
see that $\lambda_1(\lk x, T)$ can be estimated from below in terms of 
$\mu_1(\lk x)$ and the numerical invariant of $T$ which we now introduce.

\begin{Definition}\label{delta}
Let $T$ be a Hadamard space.
Suppose that collections of distinct points $\{ v_1,\dots, v_m \}\subset T$ 
and positive real numbers $\{ t_1,\dots, t_m \}$ with 
$\sum_{i=1}^m t_i = 1$ are given, and let $\overline{v}\in T$ be 
the barycenter of $\{v_i\}$ with weight $\{t_i\}$.
A {\em realization} of $\{v_i\}$ is a collection of vectors 
$\{{\bf v}_1,\dots,{\bf v}_m\} \subset \mathbb{R}^N$ for some $N$ such that
$$
||{\bf v}_i|| = d_T(\overline{v}, v_i),\quad
||{\bf v}_i-{\bf v}_j|| \leq d_T(v_i, v_j).
$$
We set
$$
\delta(\{v_i\}, \{t_i\}) = \inf
\left[\biggl|\biggl|\sum_{i=1}^m t_i{\bf v}_i\biggr|\biggr|^2 
\biggm/ \sum_{i=1}^m t_i ||{\bf v}_i||^2\right] \in [0,1],
$$
where the infimum is taken over all realizations $\{{\bf v}_i\}$
of $\{v_i\}$.
We then define
$$
\delta(T) = \sup \delta(\{v_i\}, \{t_i\}),
$$
where the supremum is taken over all collections $\{ v_i \}$, $\{ t_i \}$.
Here, if we restrict the choices of $\{ v_i \}$, $\{ t_i \}$ 
to those with barycenter at a given $v \in T$, we denote the 
corresponding number by $\delta(T, v)$.
Clearly $0\leq \delta(T, v) \leq \delta(T) \leq 1$ and
$\delta(T) = \sup_{v\in T} \delta(T, v)$.
We say that $T$ is {\em flexible} if $\delta(T) = 0$.
\end{Definition}

\begin{Lemma}\label{delta_lemma}
Let $T$ be a Hadamard space.
For $v\in T$, we have
$$
\delta(T,v) \leq \delta(TC_vT, 0_v).
$$
\end{Lemma}

\begin{proof}
Suppose that we are given distinct points $v_1,\dots,v_m \in T$ and 
positive real numbers $t_1,\dots,t_m$ with $\sum_{i=1}^m t_i = 1$
such that the barycenter of $\{ v_i \}$ with weight $\{ t_i \}$ 
coincides with $v$, and let $w_i = \pi_v(v_i)\in TC_vT$.
Then one can verify that the barycenter of $\{ w_i \}$ with weight 
$\{ t_i \}$ coincides with $0_v$ (c.f. the proof of Proposition
\ref{variation_formula}).
Let $\{ {\bf w}_i \} \subset \mathbb{R}^N$ be a realization of $\{w_i\}$.
Then since $\pi_v$ is distance non-increasing, $\{ {\bf w}_i \}$ is also
a realization of $\{v_i\}$.
Therefore, $\delta(\{v_i\}, \{t_i\}) \leq \delta(\{w_i\}, \{t_i\})$, 
and we may conclude
$\delta(T, v) \leq \delta(TC_vT, 0_v)$.
This completes the proof of Lemma \ref{delta_lemma}.
\end{proof}

\begin{Proposition}\label{wang_invariant}
Let $T$ be a Hadamard space.
Then we have
\begin{equation}\label{wang_estimate}
(1-\delta(T))\mu_1(\lk x) \leq \lambda_1(\lk x, T) \leq \mu_1(\lk x).
\end{equation}
In particular, if $T$ is flexible, 
then we have 
$\lambda_1(\lk x, T) = \mu_1(\lk x)$.
\end{Proposition}
\begin{proof}
To obtain the second inequality of \eqref{wang_estimate}, 
simply note that $T$ contains 
a geodesic, that is, an isometrically embedded line segment, 
and we can easily transform any $\mathbb{R}$-valued $0$-cochain 
to a map from $(\lk x)(0)$ into the geodesic without changing 
the Rayleigh quotient.

We shall verify the first inequality of \eqref{wang_estimate}.
Given a nonconstant map $\varphi : (\lk x)(0) \longrightarrow T$,
there exist vectors $\{{\bf v}(y) \mid y\in (\lk x)(0)\} \subset 
\mathbb{R}^N$, where $N= \#(\lk x)(0)$, such that
$$
||{\bf v}(y)|| = 
d_T(\overline{\varphi}, \varphi(y)),\quad
||{\bf v}(y) - {\bf v}(z)|| \leq d_T (\varphi(y), \varphi(z)),
$$
$$
m(x)||\overline{\bf v}||^2 \leq 
\delta(T) \sum_{y\in (\lk x)(0)} m(x,y) ||{\bf v}(y)||^2,
$$
where $\overline{\bf v} = \sum_{y\in (\lk x)(0)}(m(x,y)/m(x)){\bf v}(y)$.
We infer from the last inequality that
$$
(1-\delta(T)) ||{\bf v}||_{L^2}^2 \leq ||{\bf v}-\overline{\bf v}||_{L^2}^2.
$$
Therefore, we obtain
\begin{eqnarray*}
\lefteqn{\frac{\frac{1}{2}\sum_{(y,z)\in (\overrightarrow{\lk x})(1)} 
m(x,y,z) d_T(\varphi(y), \varphi(z))^2}{
\sum_{y\in (\lk x)(0)} m(x,y) d_T(\overline{\varphi}, 
\varphi(y))^2}}\\
&\geq& \frac{\frac{1}{2}\sum_{(y,z)\in (\overrightarrow{\lk x})(1)} 
m(x,y,z) ||{\bf v}(y) - {\bf v}(z)||^2}{
\sum_{y\in (\lk x)(0)} m(x,y) ||{\bf v}(y)||^2}\\
&=& \frac{||d{\bf v}||_{L^2}^2}{||{\bf v}||_{L^2}^2} \,\,\geq\,\, 
(1-\delta(T)) \frac{||d{\bf v}||_{L^2}^2}{
||{\bf v}-\overline{\bf v}||_{L^2}^2}.
\end{eqnarray*}
It follows that $\lambda_1(\lk x, T) \geq (1-\delta(T)) \mu_1(\lk x)$.
This completes the proof of Proposition \ref{wang_invariant}.
\end{proof}

\begin{Remark}
In Wang's fixed-point theorem mentioned in \S 3, one can replace
the assumption on $\lambda_1$ by the following:
\begin{equation}\label{general_criterion}
\mu_1(\lk x) > \frac{1}{2}
\left(1-\delta(TC_pY, 0_p)\right)^{-1}
\quad \mbox{for all $x\in X(0)$ and all $p\in Y$}.
\end{equation}
Indeed, let $f : X(0) \longrightarrow Y$ be the energy-minimizing
$\rho$-equivariant map of Theorem \ref{existence}, and set
$F_x(y) = \pi_{f(x)} (f(y))$ for $x\in X(0)$ and $y\in (\lk x)(0)$.
By Proposition \ref{variation_formula}, for any $x\in X(0)$
the barycenter of 
$\{ F_x(y) \mid y\in (\lk x)(0) \} \subset TC_{f(x)}Y$ with weight 
$\{ m(x,y)/m(x) \mid y\in (\lk x)(0) \}$ coincides with the origin 
of $TC_{f(x)}Y$.
Therefore, by the above proof, we have
\begin{eqnarray*}
\lefteqn{\sum_{(y,y') \in (\overrightarrow{\lk x})(1)}
m(x,y,y')\, d_{TC_{f(x)}Y}(F_x(y), F_x(y'))^2}\\
&\geq& 2\left(1-\delta(TC_pY, 0_p)\right)\mu_1(\lk x)
\sum_{y\in (\lk x)(0)} m(x,y)\, d_{TC_{f(x)}Y}(0_{f(x)}, F_x(y))^2.
\end{eqnarray*}
As before, this implies that $f$ should be a constant map.

The above argument shows that under the same assumption,
any harmonic $\rho$-equivariant map $f : X(0) \longrightarrow Y$ 
should be a constant map.
\end{Remark}

As a special case of Theorem \ref{fixed_point_theorem1}, we have

\begin{Theorem}\label{flexible_case}
Let $X$ be a simplicial complex equipped with an admissible weight $m$.
We assume that $\lk x$ is connected and 
\begin{equation}\label{spectral_condition}
\mu_1(\lk x) > \frac{1}{2}\quad \mbox{for all $x\in X(0)$.}
\end{equation}
Let $\Gamma$ be a finitely generated group acting by automorphisms,
properly discontinuously and cofinitely on $X$.
Let $Y$ be a Hadamard space all of whose tangent cones are flexible.
Then for any homomorphism $\rho:\Gamma \longrightarrow \isom(Y)$, 
$\rho(\Gamma)$ has a fixed point in $Y$.
\end{Theorem}

There are plenty of simplicial complexes satisfying
\eqref{spectral_condition} (with respect to the standard admissible 
weights); see Remark \ref{remark_fixed_point_theorem}.
We list some explicit examples.

\begin{Example}\label{BT_building} (Euclidean buildings)

A building is a simplicial complex which can be expressed as 
the union of subcomplexes (called apartments) satisfying 
a certain set of axioms.
In particular, it is requested that each apartment be isomorphic
to a Coxeter complex (a simplicial complex canonically 
associated to a Coxeter group) of the same type.
(The monograph \cite{Bro} is a basic reference for the theory 
of buildings.)
A building is called Euclidean if its apartments
are isomorphic to a Euclidean Coxeter complex.
A Euclidean building is contractible in general.
A well-known example of Euclidean building is given by
the one associated with the simple algebraic group 
$PGL(n, \mathbb{Q}_p)$, where $p$ is a prime and 
$\mathbb{Q}_p$ is the $p$-adic number field.
The building (more precisely, its vertex set) is the 
quotient space
$$
X = PGL(n,\mathbb{Q}_p)/PGL(n,\mathbb{Z}_p),
$$
where $\mathbb{Z}_p$ is the $p$-adic integer ring.
The dimension of $X$ is $n-1$, and the apartments are
isomorphic to a Euclidean Coxeter complex of type
$\widetilde{A}_{n-1}$.
If $n=2$, $X$ is a regular tree of degree $p+1$.
If $n=3$, $X$ is two-dimensional and the links of vertices 
are all isomorphic to the same graph $\mathcal{G}$; 
it is a regular bipartite graph of degree $p+1$ with 
$2(p^2+p+1)$ vertices and $(p+1)(p^2+p+1)$ edges.
If $p=2$, $\mathcal{G}$ is as in the following picture:

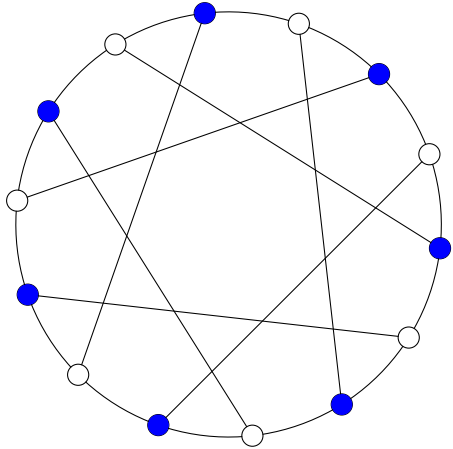
\vspace{.5cm}
\hspace{4cm}
\includegraphics{link.1}
\vspace{.5cm}

\noindent
By a result of Feit-Higman \cite{FeiHig}, the first nonzero
eigenvalue of the Laplacian of $\mathcal{G}$ is
$1-\sqrt{p}/(p+1) > 1/2$.
Examples of $\Gamma$ are supplied by cocompact lattices 
in $PGL(3, \mathbb{Q}_p)$.
\end{Example}

\begin{Example}\label{BS_complex} (Ballmann-\'Swi\c{a}tkowski complexes
\cite[Theorem 2]{BalSwi})

Let $H$ be a finite group, $S\subset H\setminus \{e\}$ a set
of generators of $H$ and $L=C(H,S)$ the Cayley graph of $H$
with respect to $S$.
Assume that the girth of $L$ (i.e., the minimal number of edges
in closed circuits of $L$) is at least $6$.
Then there exists a contractible two-dimensional simplicial
complex $X$ such that the links of all vertices of $X$ are
isomorphic to $L$.
Moreover, if $\langle S,R \rangle$ is a presentation of $H$,
then the group $\Gamma$ given by the presentation 
$$
\langle S \cup \{\tau\} \mid R \cup \{\tau^2\} \cup 
\{(s\tau)^3 \mid s\in S\} \rangle
$$
acts by automorphisms, properly discontinuously and cofinitely on $X$.

Sarnak \cite[Chapter 3]{Sar} describes some explicit examples of 
Ramanujan graphs as Cayley graphs of finite groups.
For most of them, the first nonzero eigenvalue of the Laplacian 
is greater than $1/2$.
Therefore, taking such groups for $H$ gives rise to examples of 
discrete groups in demand.
\end{Example}

A tangent cone is flexible if it is isometric to a Hilbert space.
On the other hand, there certainly exist Hadamard spaces with 
non-Hilbertian, flexible tangent cones.

\begin{Example}\label{tree} (trees)\quad
Let $Y$ be a locally finite tree in which each vertex belongs to
at least three edges.
By defining the length of each edge to be one, for example, 
$Y$ becomes a Hadamard space.
If $p$ is an interior point of an edge of $Y$, the tangent cone
at $p$ is isometric to a line.
If $p$ is a vertex of $Y$, then the tangent cone at $p$ is
isometric to an $n$-pod $T_n$ for some $n$, which is the union of 
$n$ half-lines with all end points identified.
The angle between vectors of $T_n$ is given by
$$
\angle(v,w) = \left\{ \begin{array}{cl}
0 & \mbox{if $v$, $w$ belong to the same half-line},\\
\pi & \mbox{otherwise}. \end{array} \right.
$$
We claim that $T_n$ is flexible.
For the sake of simplicity, we shall prove this for $n=3$.
The argument easily extends to general $n$.

We denote the half-lines of $T_3$ by $H_s$, $s=1,2,3$, and
the origin of $T_3$ by $0$.
By Lemma \ref{delta_lemma}, it suffices to show that 
$\delta(T_3, 0) = 0$.
Suppose that we are given nonzero vectors $v_1,\dots, v_m \in 
T_3$ and positive real numbers $t_1,\dots, t_m$ with
$\sum_{i=1}^m t_i = 1$ such that 
the barycenter of $\{v_i\}$ with weight $\{t_i\}$ coincides
with $0$.
One can verify, by computation similar to the one in the proof
of Proposition \ref{variation_formula}, that the inequality 
$\sum_{i=1}^m t_i \langle v_i, w \rangle \leq 0$ holds for all 
$w\in T_3$.
Let $I_s = \{ i \mid v_i \in H_s \}$ and set
$A_s = \sum_{i\in I_s} t_i |v_i|$ for $s=1,2,3$.
Then the last condition is translated to the system of triangle 
inequalities
$$
A_1\leq A_2 + A_3,\quad A_2\leq A_3 + A_1,\quad
A_3\leq A_1 + A_2.
$$ 
But this means that there exists a (possibly degenerate) triangle 
in $\mathbb{R}^2$ whose sides have length $A_s$.
In other words, there exist unit vectors ${\bf e}_s$ such that 
$\sum_{s=1}^3 A_s {\bf e}_s = {\bf 0}$.
We now set ${\bf v}_i = |v_i| {\bf e}_s$ if $i\in I_s$.
Then $\{{\bf v}_i\}$ is a realization of $\{v_i\}$ and satisfies
$\sum_{i=1}^m t_i {\bf v}_i = {\bf 0}$.
Thus $\delta(T_3, 0) = 0$, and $T_3$ is flexible.

Since Definition \ref{delta} concerns only a finite number of points,
the above claim remains the case if $n$ is infinite.
Therefore, the tangent cones of an $\mathbb{R}$-tree are also all 
flexible.
\end{Example}

\begin{Example}\label{product_of_trees} (product of trees)\quad
Let $Y_\nu$, $\nu=1,\dots,k$, be trees as in Example \ref{tree}.
Then the tangent cones of the product $\prod_{\nu=1}^k Y_\nu$
are all flexible.
This is a consequence of the following general fact.
\end{Example}

\begin{Proposition}\label{flexibility_of_product}
Let $T_1$ and $T_2$ be two Hadamard spaces.
Then the following inequality holds{\rm :}
\begin{equation}\label{delta_of_product}
\delta(T_1 \times T_2) \leq \max \{\delta(T_1), \delta(T_2) \}.
\end{equation}
In particular, if $T_\nu$ is flexible for each $\nu$, then 
$T_1 \times T_2$ is also flexible.
\end{Proposition}
\begin{proof}
Suppose that we are given distinct points $v_1,\dots,v_m \in
T_1\times T_2$ and positive real numbers $t_1,\dots,t_m$ with 
$\sum_{i=1}^m t_i = 1$, and let $\overline{v}$ be the barycenter
of $\{v_i\}$ with weight $\{t_i\}$.
Write $v_i = (v_i^{(1)}, v_i^{(2)})$ and $\overline{v}
= (\overline{v}^{(1)}, \overline{v}^{(2)})$.
It is easy to see that the barycenter of $\{v_i^{(\nu)}\}$
with weight $\{t_i\}$ coincides with $\overline{v}^{(\nu)}$
for each $\nu$.
By the definition of $\delta$, for each $\nu$ there exists a realization
$\{ {\bf v}_i^{(\nu)}\} \subset \mathbb{R}^{N_\nu}$ of
$\{v_i^{(\nu)}\}$ satisfying
$$
\biggl|\biggl|\sum_{i=1}^m t_i{\bf v}_i^{(\nu)}\biggr|\biggr|^2
\leq \delta(T_\nu) \sum_{i=1}^m t_i ||{\bf v}_i^{(\nu)}||^2.
$$
Let ${\bf v}_i = ({\bf v}_i^{(1)}, {\bf v}_i^{(2)})$.
Then it is easy to see that 
$\{ {\bf v}_i \} \subset \mathbb{R}^{N_1+N_2}$ is a realization 
of $\{ v_i \}$ 
and 
$$
\biggl|\biggl|\sum_{i=1}^m t_i{\bf v}_i\biggr|\biggr|^2
\leq \max_\nu \delta(T_\nu) \sum_{i=1}^m t_i ||{\bf v}_i||^2.
$$
Thus \eqref{delta_of_product} follows, and the proof of Proposition 
\ref{flexibility_of_product} is completed.
\end{proof}

Applying Theorem \ref{flexible_case} to the examples above, 
we obtain
\begin{Corollary}\label{fixed_point_theorem2}
Let $X$ and $\Gamma$ be as in Theorem \ref{flexible_case}.
Let $Y$ be one of the following spaces{\rm :}
\begin{enumerate}
\renewcommand{\theenumi}{\roman{enumi}}
\renewcommand{\labelenumi}{(\theenumi)}
\item a Hadamard space all of whose tangent cones are isometric 
to closed convex cones of Hilbert spaces,
\item an $\R$-tree,
\item a product of Hadamard spaces which are either of type {\rm (i)} 
or {\rm (ii)}.
\end{enumerate}
Then for any homomorphism $\rho:\Gamma \longrightarrow \isom(Y)$, 
$\rho(\Gamma)$ has a fixed point in $Y$.
\end{Corollary}

\begin{Remark}\label{remark_fixed_point_theorem}
When $Y$ is a (finite-dimensional) Hadamard manifold and $\rho$ is reductive, 
the assertion of the theorem was proved by Wang \cite{Wan1}.

 Taking $Y$ to be a Hilbert space in Theorem \ref{fixed_point_theorem2}, 
 we see that any isometric action of $\Gamma$ on a Hilbert space has a 
 fixed point. 
 This implies that $\Gamma$ has Kazhdan's Property (T)
 (see \cite{HarVal}). 
 Ballmann and \'Swi\c{a}tkowski \cite[Corollary 1]{BalSwi} gave a proof
 of this fact, assuming that $X$ is contractible.
 Our argument shows that the connectedness of $\lk x$ is sufficient.

  Our result also recovers a part of \.Zuk's result \cite[Theorem 1]{Zuk2},
 since one can construct a two-dimensional simplicial complex on which
 $\Gamma$ acts by automorphisms, properly discontinuously and cofinitely,
 so that each $\lk x$ is isomorphic to \.Zuk's $L(S)$, where $S$ is a
 generator set for $\Gamma$. 
 Together with \.Zuk's result \cite[Theorem 4]{Zuk2}, this argument assures
 plenty of simplicial complexes satisfying \eqref{spectral_condition}.

It is known that any isometric action of a locally 
compact group having Kazhdan's property (T) on an $\R$-tree has a fixed point 
(see \cite{HarVal}).
 On the other hand, our theorem claims that any isometric action of $\Gamma$,
 for example, on a (possibly infinite-dimensional) Hadamard manifold must have 
 a fixed point, and this is not true of a general group satisfying Kazhdan's 
 property (T). (Consider an irreducible lattice of a real semisimple 
 Lie group of rank $\geq 2$.)
 
Gromov \cite[p.~16]{Gro1} asked whether a finitely generated group
of finite cohomological dimension is always a fundamental group
of a complete Riemannian manifold of nonpositive sectional curvature.
Our result gives a negative answer to his question.
For example, a cocompact lattice in $PGL(n, \mathbb{Q}_p)$ has a torsion-free
subgroup $\Gamma$ acting very freely on the associated Euclidean building.
This $\Gamma$ is finitely generated and has finite cohomological dimension.
Theorem \ref{fixed_point_theorem2} implies, however, that it cannot act 
properly discontinuously on a Hadamard manifold.

As mentioned in the introduction, Gromov \cite{Gro2} also proves results similar
to Theorem \ref{fixed_point_theorem1} and Corollary \ref{fixed_point_theorem2}.
\end{Remark}

We shall now give an example of a Hadamard space having 
non-flexible tangent cones.

\begin{Example}\label{A_2_tilde_building}
We shall take up the $\widetilde{A}_2$-building of Example
\ref{BT_building} again.
In general, a Euclidean building can be equipped with a 
distance by transplanting the Euclidean distance onto
each apartment, and the building becomes a Hadamard space
with this distance (see \cite[Chapter 6]{Bro}).

Let $Y$ be the $\widetilde{A}_2$-building of Example
\ref{BT_building}:
$$
Y  = PGL(3,\mathbb{Q}_r)/PGL(3,\mathbb{Z}_r).
$$
If $p\in Y$ is not a vertex, then $TC_pY$ is flexible, 
as it is isometric to
either a Euclidean plane or a product of ($r+1$)-pod with a line.
If $p\in Y$ is a vertex, $TC_pY$ is isometric to the metric cone 
$C(\mathcal{G})$
over the graph $\mathcal{G}$ which we equip with a metric by assigning 
length $\pi/3$ to each edge.
(The distance of $C(\mathcal{G})$ is defined as in Definition 
\ref{tangent cone} (3).) 
We observe
\begin{equation}\label{lower_bound_for_delta}
\delta(TC_pY, 0_p) \geq 1-\frac{1}{2(1-\sqrt{r}/(r+1))} 
= \frac{(\sqrt{r}-1)^2}{2(r-\sqrt{r}+1)}.
\end{equation}
In particular, $TC_pY$ is not flexible.
To see this, let $X=Y$, $\Gamma$ 
a cocompact lattice of $PGL(3,\mathbb{Q}_p)$ and $\rho$ the 
composition of the inclusion maps $\Gamma \hookrightarrow 
PGL(3,\mathbb{Q}_p) \hookrightarrow \isom(Y)$.
Then it is easy to verify that the identity map $\id : X \longrightarrow Y$ 
is a harmonic $\rho$-equivariant map.
This means that the inequality \eqref{general_criterion} cannot hold.
Since $\mu_1(\lk x) = 1-\sqrt{r}/(r+1)$ for all $x\in X(0)$, we obtain 
\eqref{lower_bound_for_delta}.
Note that $\delta(TC_pY, v) = 0$ if $v\neq 0_p$.
\end{Example}

We shall make some more observation concerning this example 
and prove a fixed-point theorem in the next section.

\section{Optimal embedding of $C(\mathcal{G})$}

Throughout this section, $Y$ denotes the $\widetilde{A}_2$-building 
of Example \ref{A_2_tilde_building} and $p$ is a vertex of $Y$.
In this section, we shall construct a certain embedding of 
the tangent cone $TC_pY$ into a Euclidean space.
We shall then use it to estimate $\delta(TC_pY) = \delta(TC_pY, 0_p)$ 
from above when $r=2$.

Recall that $TC_pY$ is isometric to the metric cone $C(\mathcal{G})$ 
over the metric graph $\mathcal{G}$. 
We number the vertices of $\mathcal{G}$ by $s = 1,\dots,2(r^2+r+1)$.
When we emphasize that they are points of $C(\mathcal{G})$, we shall denote
them by $e_s$.
Note that the barycenter of $\{ e_s \}$ coincides with the origin of 
$C(\mathcal{G})$.
We denote the combinatorial distance on $\mathcal{G}$ by $d$.
Then since $\mathcal{G}$ is a $(r+1)$-regular, bipartite graph with
girth $6$, there are $(r+1)(r^2+r+1)$ (resp. $(r^2+r)(r^2+r+1)$, 
$r^2(r^2+r+1)$) pairs of vertices with $d=1$ (resp. $d=2$, $d=3$).
Note that two vertices $s,t$ are of the same type if and only if
$d(s,t) =2$.
By a {\em chamber}, we mean a subset of $C(\mathcal{G})$ 
which is the cone over an edge of $\mathcal{G}$.

We shall construct a map $\iota$ from $C(\mathcal{G})$ into the
Euclidean space $\R^N$, where $N = r^2+r+1$, such that 
\begin{enumerate}
\renewcommand{\theenumi}{\roman{enumi}}
\renewcommand{\labelenumi}{(\theenumi)}
\item $\iota$ satisfies 
$$
||\iota(v)|| = d_{C(\mathcal{G})}(0, v),\quad
||\iota(v)-\iota(w)|| \leq d_{C(\mathcal{G})}(v, w)
$$
for all $v,w\in C(\mathcal{G})$, where $0$ is the origin of 
$C(\mathcal{G})$;
\item the barycenter of $\{ \iota(e_s) \}$ is as close to the origin 
${\bf 0}\in \R^N$
as the one associated with any other map from $C(\mathcal{G})$
into a Euclidean space of any dimension, satisfying the condition (i).
\end{enumerate}
In fact, the map $\iota$ constructed below has the additional property 
that it is distance preserving when restricted to each chamber.

To find such an $\iota$, we shall proceed as follows.
Let $V = \oplus_s \R {\bf e}_s$ be the vector space 
with basis $\{ {\bf e}_s \}$. 
We introduce a family of symmetric bilinear forms 
$\{ \langle\cdot,\cdot\rangle_{a,b} \mid -1\leq a,b\leq 1 \}$ on $V$ 
by $\langle {\bf e}_s, {\bf e}_s \rangle_{a,b} = 1$ and
$$
\langle {\bf e}_s, {\bf e}_t \rangle_{a,b} = \left\{\begin{array}{ccc}
1/2 &\mbox{if}& d(s,t)=1,\\
a &\mbox{if}& d(s,t)=2,\\
b &\mbox{if}& d(s,t)=3.
\end{array} \right.
$$
Then 
\begin{equation}\label{barycenter_size}
\biggl\langle \sum_s {\bf e}_s, \sum_s {\bf e}_s \biggr\rangle_{a,b}
= (r+3)(r^2+r+1) + 2r(r^2+r+1)\{(r+1)a+rb\}.
\end{equation}
We shall minimize this quantity under the constraint that 
$\langle\cdot,\cdot\rangle_{a,b}$ is positive semidefinite.
The Gram matrix $G_{a,b} = (\langle {\bf e}_s,{\bf e}_t \rangle_{a,b})$ 
is related to the adjacency matrix $A$ of $\mathcal{G}$ by 
$$
G_{a,b} = \left(\begin{array}{cc} aJ & bJ\\ bJ & aJ \end{array}\right)
+ (1-a) I + (1/2 - b) A,
$$
where $J$ is the matrix of size $r^2+r+1$ all of whose entries
are $1$.
Note that the eigenvalues of the first matrix on the right-hand side are
$(r^2+r+1)(a\pm b)$ with eigenvectors
${\bf v}^{\pm} = \mbox{}^t(1,\dots,1,\pm 1,\dots,\pm 1)$ 
and $0$ with multiplicity $2(r^2+r)$.
Since the eigenvalues of $A$ are $\pm(r+1)$ with eigenvectors 
${\bf v}^{\pm}$
and $\pm \sqrt{r}$ with multiplicities $r^2+r$ (see \cite{FeiHig, Zuk2}), 
those of $G_{a,b}$ are given by
$$
(r^2+r+1)(a\pm b) + (1-a) \pm (1/2-b)(r+1)
$$
with multiplicities $1$ and
$$
(1-a)\pm (1/2-b)\sqrt{r}
$$
with  multiplicities $r^2+r$.
Under the constraint that these are nonnegative, the quantity
\eqref{barycenter_size} takes minimum value
$2(r^2+r+1)(r^2 + 1- (r+1)\sqrt{r})$ at
$$
(a,b) = \left(\frac{r-1-\sqrt{r}}{2r}, \frac{r^2-r-(r+1)\sqrt{r}}{2r^2}
\right).
$$
With this $(a,b)$, $G_{a,b}$ has zero eigenvalue with multiplicity $r^2+r+1$
and positive eigenvalues $r^2 + 1 - (r+1)\sqrt{r}$, $(r+1+\sqrt{r})/r$
with multiplicities $1$, $r^2+r$ respectively.
We denote the sum of the eigenspaces corresponding to these two positive
eigenvalues of $G_{a,b}$ by $W$, and denote the projection from $V$ onto
$W$ by $\pi$.
By restricting $\langle\cdot,\cdot\rangle_{a,b}$ onto $W$, we obtain 
a Euclidean space of dimension $r^2+r+1$, and $\pi$ preserves
the inner products.
Let $\iota$ be the composition of the maps $C(\mathcal{G}) \longrightarrow
V \longrightarrow W$, where the first one is the natural inclusion
sending $e_s$ to ${\bf e}_s$.
Then $\iota$ clearly satisfies the condition (i) above, and
setting ${\bf e}_s' = \iota(e_s) = \pi({\bf e}_s)$, we have
$$
\biggl|\biggl| \sum_s {\bf e}_s' \biggr|\biggr|^2
= 2(r^2+r+1)(r^2 + 1- (r+1)\sqrt{r}),
$$
or
$$
\biggl|\biggl| \sum_s t_s {\bf e}_s' \biggr|\biggr|^2 
= \frac{(\sqrt{r}-1)^2}{2(r-\sqrt{r}+1)} \sum_s t_s ||{\bf e}_s'||^2,
$$
where $t_s = 1/2(r^2+r+1)$ for all $s$.
Together with the argument in Example \ref{A_2_tilde_building},
this implies
$$
\delta(\{e_s\}, \{t_s\}) = \frac{(\sqrt{r}-1)^2}{2(r-\sqrt{r}+1)}.
$$
Thus $\iota$ satisfies the condition (ii) above also.

The maximal eigenvalue of $G_{a,b}$ is given by
$$
R(G_{a,b}) = \left\{\begin{array}{ll} (3+\sqrt{2})/2, & \mbox{$r=2$},\\
r^2 + 1 - (r+1)\sqrt{r}, & \mbox{otherwise}.\end{array}\right.
$$
Given any collection $\{ A_s \}$ of positive real numbers, we have
\begin{eqnarray*}
\biggl|\biggl| \sum_s A_s {\bf e}_s' \biggr|\biggr|^2
&=& (A_1,\dots) G \,\,\mbox{}^t (A_1,\dots)\\
&\leq& R(G) \sum_s {A_s}^2.
\end{eqnarray*}
Therefore,
$$
\biggl|\biggl| \sum_s t_s A_s {\bf e}_s' 
\biggr|\biggr|^2
\leq \frac{R(G)}{2(r^2+r+1)} \sum_s t_s ||A_s {\bf e}_s'||^2,
$$
where $t_s$ are as above.
This means that
$$
\delta(\{A_s e_s\}, \{t_s\}) \leq
\frac{(\sqrt{r}-1)^2}{2(r-\sqrt{r}+1)}
$$
unless $r=2$, which partially confirms the following

\begin{Conjecture}\label{conjecture}
$$
\delta(TC_pY, 0_p) = \frac{(\sqrt{r}-1)^2}{2(r-\sqrt{r}+1)}.
$$
\end{Conjecture}

This conjecture must be left unanswered as a future problem. 
On the other hand, when $r=2$, we can use the above map $\iota$ 
to give a meaningful upper estimate of 
$\delta(TC_pY) = \delta(TC_pY, 0_p)$. 
(The argument below gives a similar estimate for larger $r$, 
but the result is not good enough in view of the application.)
Set $G = G_{a,b}$ with the values of $a,b$ as above and 
$\overline{G} = G_{-1/2,-1}$.
The latter is the Gram matrix associated with the symmetric bilinear 
form $\langle\langle\cdot,
\cdot\rangle\rangle = \langle\cdot,\cdot\rangle_{-1/2,-1}$ on $V$. 
Recall that the adjacency matrix $A$ of the graph $\mathcal{G}$
has two one-dimensional eigenspaces
$W^{(1)} = \R\, {\bf v}^+$
and
$Z^{(1)} = \R\, {\bf v}^-$,
and two six-dimensional eigenspaces $W^{(6)}$ and $Z^{(6)}$.
Here $Z^{(6)}$ is chosen so that $Z^{(1)}\oplus Z^{(6)}$
is the zero eigenspace of $G$.
Thus $G$ has two positive eigenvalues $5-3\sqrt{2}$ and 
$(3+\sqrt{2})/2$ with the corresponding eigenspaces $W^{(1)}$ 
and $W^{(6)}$ respectively.
On the other hand, $\overline{G}$ has eigenvalues 
$-9/2$, $3(1+\sqrt{2})/2$, $1/2$ and $3(1-\sqrt{2})/2$ 
with the corresponding eigenspaces $W^{(1)}$, $W^{(6)}$, 
$Z^{(1)}$ and $Z^{(6)}$ respectively. 

Now suppose that nonzero vectors $\{ v_1,\dots v_m \} \subset
C(\mathcal{G})$ and positive real numbers $\{ t_1,\dots,t_m \}$ 
are given so that the 
barycenter of $\{ v_i \}$ with weight $\{ t_i \}$ coincides with
the origin of $C(\mathcal{G})$.
Set ${\bf v} = \sum_{i=1}^m t_i {\bf v}_i \in V$.
Hereafter we shall identify $\sum_s a_s {\bf e}_s \in V$ with 
the corresponding column vector $\mbox{}^t (a_1,\dots,a_{14})$.
It is easy to see that the barycentric condition implies
$\overline{G} {\bf v}$ is a negative vector, that is, the components
of $\overline{G} {\bf v}$ are all nonpositive.
In particular, we have $\langle\langle {\bf v}, {\bf v} \rangle\rangle 
= \mbox{}^t{\bf v} \overline{G} {\bf v} \leq 0$
since ${\bf v}$ is a positive vector.
Let ${\bf v} = {\bf w}^{(1)} + {\bf w}^{(6)} + {\bf z}^{(1)} 
+ {\bf z}^{(6)}$ and ${\bf v}_i = {\bf w}^{(1)}_i + {\bf w}^{(6)}_i 
+ {\bf z}^{(1)}_i + {\bf z}^{(6)}_i$, $i=1,\dots,m$, be the 
decompositions corresponding to the decomposition 
$V = W^{(1)} \oplus W^{(6)} \oplus Z^{(1)} \oplus Z^{(6)}$.
Since 
$$
\langle\langle {\bf v}, {\bf v} \rangle\rangle 
= \langle\langle {\bf w}^{(1)}, {\bf w}^{(1)} \rangle\rangle 
+ \langle\langle {\bf w}^{(6)}, {\bf w}^{(6)} \rangle\rangle 
+ \langle\langle {\bf z}^{(1)}, {\bf z}^{(1)} \rangle\rangle 
+ \langle\langle {\bf z}^{(6)}, {\bf z}^{(6)} \rangle\rangle 
\leq 0
$$
and 
$\langle\langle {\bf z}^{(1)}, {\bf z}^{(1)} \rangle\rangle \geq 0$,
we have
$$
\langle\langle {\bf w}^{(6)}, {\bf w}^{(6)} \rangle\rangle
\leq - \langle\langle {\bf w}^{(1)}, {\bf w}^{(1)} \rangle\rangle
- \langle\langle {\bf z}^{(6)}, {\bf z}^{(6)} \rangle\rangle.
$$
Using this, we estimate
\begin{eqnarray*}
\langle {\bf v}, {\bf v} \rangle 
&=& \langle {\bf w}^{(1)}, {\bf w}^{(1)} \rangle
+  \langle {\bf w}^{(6)}, {\bf w}^{(6)} \rangle \\
&=& \frac{5-3\sqrt{2}}{-9/2} 
\langle\langle {\bf w}^{(1)}, {\bf w}^{(1)} \rangle\rangle 
+ \frac{3+\sqrt{2}}{3(1+\sqrt{2})}
\langle\langle {\bf w}^{(6)}, {\bf w}^{(6)} \rangle\rangle \\
&\leq& \frac{2\sqrt{2}-1}{3}\left[ \frac{2\sqrt{2}+1}{3}
(- \langle\langle {\bf w}^{(1)}, {\bf w}^{(1)} \rangle\rangle)
+ (- \langle\langle {\bf z}^{(6)}, {\bf z}^{(6)} \rangle\rangle)
\right] \\
&\leq& \frac{2\sqrt{2}-1}{3}\sum_{i=1}^m t_i \left[
\frac{2\sqrt{2}+1}{3}
(- \langle\langle {\bf w}_i^{(1)}, {\bf w}_i^{(1)} \rangle\rangle)
+ (- \langle\langle {\bf z}_i^{(6)}, {\bf z}_i^{(6)} \rangle\rangle)
\right].
\end{eqnarray*}
For the last inequality, we have used the Cauchy-Schwarz inequality.

It remains to estimate 
$-\langle\langle {\bf w}_i^{(1)}, {\bf w}_i^{(1)} \rangle\rangle$
and $-\langle\langle {\bf z}_i^{(6)}, {\bf z}_i^{(6)} \rangle\rangle$
from above in terms of 
$\langle\langle {\bf v}_i, {\bf v}_i \rangle\rangle
= \langle {\bf v}_i, {\bf v}_i \rangle$.
Note that each ${\bf v}_i$ has the form $a{\bf e}_s + b{\bf e}_t$,
where $a,b\geq 0$ and $s$, $t$ are neighboring vertices of 
$\mathcal{G}$.
We assume ${\bf v}_i = {\bf e}_s$ for the moment and express
${\bf w}_i^{(1)}$ and 
${\bf z}_i^{(6)}$ as linear combinations of ${\bf e}_u$'s.
Computation of ${\bf w}_i^{(1)}$ is easy;
we obtain
$$
{\bf w}_i^{(1)} = \frac{1}{14} \sum_u {\bf e}_u.
$$
To compute ${\bf z}_i^{(6)}$, we observe that the coefficients 
of ${\bf e}_u$ depend only on $d(s,u)$.
Indeed, the group $H= GL(3, \Z/2\Z)$ acts on $\mathcal{G}$
by automorphisms, as $\mathcal{G}$ is the incidence graph
of the finite projective plane $\mathbb{P}^2(\Z/2\Z)$.
Therefore, the  induced action of $H$ on $V$ leaves $Z^{(6)}$ 
invariant, and
the projection $V \longrightarrow Z^{(6)}$ commutes with the 
$H$-actions.
It follows that ${\bf z}_i^{(6)}$ is invariant by the stabilizer
$H_s$ of $s$.
It is easy to see from the incidence-graph picture of $\mathcal{G}$
that $H_s$ acts on the set $\{ u \in \mathcal{G}(0) \mid d(s,u)=d \}$ 
transitively for each $d\in \{1,2,3\}$ (see \cite[p.~658]{Zuk2}).
The claim now follows.
By using this fact, we can conclude
$$
{\bf z}_i^{(6)} = \frac{3}{7} {\bf e}_s 
+ \frac{\sqrt{2}}{7} \sum_{d(s,u)=1} {\bf e}_u 
- \frac{1}{14} \sum_{d(s,u')=2} {\bf e}_{u'}
- \frac{3\sqrt{2}}{28} \sum_{d(s,u'')=3} {\bf e}_{u''}.
$$
It follows that for ${\bf v}_i = a{\bf e}_s + b{\bf e}_t$, 
we have 
$\langle\langle {\bf v}_i, {\bf v}_i \rangle\rangle 
= \langle {\bf v}_i, {\bf v}_i \rangle = a^2 + ab + b^2$ 
and
$$
\langle\langle {\bf w}_i^{(1)}, {\bf w}_i^{(1)} \rangle\rangle
= - \frac{9}{28} (a+b)^2,\quad
\langle\langle {\bf z}_i^{(6)}, {\bf z}_i^{(6)} \rangle\rangle
= - \frac{3\sqrt{2}-3}{14} (3 a^2 - 2\sqrt{2} ab + 3 b^2).
$$
Therefore, we conclude
\begin{eqnarray*}
\langle {\bf v}, {\bf v} \rangle
&\leq& \frac{2\sqrt{2}-1}{3}\sum_{i=1}^m t_i \left[
\frac{2\sqrt{2}+1}{3} \frac{9}{28}(a+b)^2 
+ \frac{3\sqrt{2}-3}{14} (3 a^2 - 2\sqrt{2} ab + 3 b^2) \right] \\
&\leq& \frac{37-18\sqrt{2}}{28} \sum_{i=1}^m t_i
\langle {\bf v}_i, {\bf v}_i \rangle.
\end{eqnarray*}
Thus 
$$
\delta(TC_pY) = \delta(TC_pY, 0_p) \leq \frac{37-18\sqrt{2}}{28} 
= 0.4122\cdots.
$$
Recall that $TC_pY$ is flexible 
if $p\in Y$ is not a vertex.
Therefore, $\delta(TC_pY) \leq (37-18\sqrt{2})/28$ for all $p\in Y$.
By Proposition \ref{wang_invariant} and Theorem \ref{fixed_point_theorem1},
we obtain

\begin{Corollary}\label{fixed_point_theorem_for_A_2_tilde_building}
Let $X$ be a simplicial complex equipped with an admissible weight $m$,
and assume that $\lk x$ is connected for all $x\in X(0)$.
Let $\Gamma$ be a finitely generated group acting by automorphisms,
properly discontinuously and cofinitely on $X$.
Let $Y$ be the $\widetilde{A}_2$-building as above{\rm :}
$$
Y  = PGL(3,\mathbb{Q}_2)/PGL(3,\mathbb{Z}_2).
$$
Suppose that 
\begin{equation}\label{spectral_condition2}
\mu_1(\lk x) > \frac{1}{2\{1-(37-18\sqrt{2})/28\}} = 0.8507\cdots
\end{equation}
for all $x\in X(0)$.
Then for any homomorphism $\rho:\Gamma \longrightarrow \isom(Y)$, 
$\rho(\Gamma)$ has a fixed point in $Y$.
\end{Corollary}

\begin{Remark}
There are many simplicial complexes satisfying \eqref{spectral_condition2};
$\widetilde{A}_2$-buildings of Example \ref{BT_building} with large $p$ are 
basic ones.
In fact, \eqref{spectral_condition2} is satisfied if the prime $p$ is 
not less than $43$.
\end{Remark}

\section{Appendix. Bochner-Weitzenb\"ock-type formula for equivariant
cochains} 

Let $X$ be a simplicial complex equipped with an admissible weight $m$,
and $\Gamma$ a finitely generated group acting by automorphisms, 
properly discontinuously and cofinitely on $X$. 
Let $\mathcal{H}$ be a Hilbert space with inner product 
$\langle \cdot, \cdot \rangle$, and $\mathcal{U}(\mathcal{H})$ the group
of unitary operators on $\mathcal{H}$. 
Consider a unitary representation 
$\rho: \Gamma \longrightarrow \mathcal{U}(\mathcal{H})$ of $\Gamma$. 
We say that a simplicial $k$-cochain $\alpha$ on $X$ with values in 
$\mathcal{H}$ is $\rho$-\textit{equivariant} if it satisfies 
$\alpha (\gamma s)= \rho(\gamma)\alpha (s)$
for all $\gamma \in \Gamma$ and $s\in \ordered{k}$, and denote 
by $C^k_{\rho}(X)$ the set of $\rho$-equivariant $k$-cochains on $X$ 
with values in $\mathcal{H}$. 

We define an inner product on $C^k_{\rho}(X)$ by 
\begin{equation*}
 \langle \alpha, \beta \rangle_{L^2} 
  =  \frac{1}{(k+1)!} \sum_{s \in \orderedf{k}} 
  \frac{m(s)}{|\Gamma_s|}\langle \alpha(s), \beta(s)\rangle, 
 \quad \alpha, \beta \in C^k_{\rho}(X).
\end{equation*}
Note that the right-hand side is independent of the choice of 
$\orderedf{k}$. 

Let $d:C^{k-1}_{\rho}(X) \longrightarrow C^k_{\rho}(X)$, $k\geq 1$, 
be the simplicial coboundary operator given by 
\begin{equation*}
d\alpha (s) =\sum_{i=0}^{k}(-1)^i \alpha(s_{(i)}),\quad 
s\in \ordered{k},
\end{equation*}
where $s_{(i)}= (x_0, \dots, \hat x_i, \dots, x_{k})$
if $s=(x_0, \dots, x_{k})\in \ordered{k}$.  

\begin{Proposition}\label{adjoint}
The adjoint $\delta:C^k_{\rho}(X) \longrightarrow C^{k-1}_{\rho}(X)$ 
of $d$ takes the form of
\begin{equation}\label{adjoint_formula}
 \delta \alpha (s)= (-1)^{k} \sum_{t \in \ordered{k}_s} 
 \frac{m(t)}{m(s)}\alpha (t), \quad s \in \ordered{k-1}, 
\end{equation}
and both $d$ and $\delta$ are bounded operators. 
\end{Proposition}
\begin{proof} 
Let $\tau_{i,j}: \ordered{k}\longrightarrow \ordered{k}$ 
be the map defined by
\begin{equation*}
 \tau_{i,j}\big((x_0, \dots, x_i,\dots, x_j,\dots, x_{k})\big) 
  = (x_0, \dots,\overset{i}{\check{x_j}},\dots, 
    \overset{j}{\check{x_i}},\dots x_{k}). 
\end{equation*}
Note that   
$\alpha (\tau_{i,j}(t))= -\alpha(t)$ for $\alpha \in C^k_{\rho}(X)$ 
if $i\not= j$
and $(\tau_{i,k}(t))_{(k)}=\tau_{i,i+1}\circ \dots \circ 
\tau_{k-2,k-1}(t_{(i)})$. 
Thus, for $\alpha \in C^{k}_{\rho}(X)$ and $\beta \in C^{k-1}_{\rho}(X)$, 
we have
\begin{eqnarray*}
\lefteqn{\langle \alpha, d\beta \rangle_{L^2}} \\
 &=& \frac{1}{(k+1)!} 
     \sum_{t \in \orderedf{k}}\frac{m(t)}{|\Gamma_t|} 
     \left\langle \alpha(t), \sum_{i=0}^{k}(-1)^i \beta(t_{(i)}) 
     \right\rangle \\
 &=& \frac{1}{(k+1)!} 
     \sum_{t \in \orderedf{k}}\frac{m(t)}{|\Gamma_t|} 
     \Bigg[\sum_{i=0}^{k-1}
     \left\langle -\alpha\left(\tau_{i,k}(t) \right), 
     (-1)^{k-1} \beta (\tau_{i,k}(t)_{(k)})\right\rangle \\
 & & \phantom{\sum_{t \in \orderedf{k}}\frac{m(t)}{|\Gamma_t|}}
   + \left\langle \alpha(t),(-1)^{k}\beta(t_{(k)})
     \right\rangle \Bigg] \\
 &=& \frac{(-1)^k}{(k+1)!} \left[\sum_{i=0}^{k-1}
     \sum_{t' \in \tau_{i,k}(\orderedf{k})}
     \frac{m(t')}{|\Gamma_{t'}|} 
     \big\langle \alpha(t'),\beta (t'_{(k)})\big\rangle 
     + \sum_{t \in \orderedf{k}}\frac{m(t)}{|\Gamma_t|}
     \big\langle\alpha(t), \beta(t_{(k)})\big\rangle \right].
\end{eqnarray*}
Note that $\tau_{i,j}(\orderedf{k})$ is also a representative set, 
and the sums over representative sets in the computation above are
independent of the choice of the representative sets. 
Therefore, we obtain 
\begin{equation}\label{rewrite}
\langle \alpha, d\beta \rangle_{L^2}  
= \frac{(-1)^{k}}{k!} \sum_{t \in \orderedf{k}}\frac{m(t)}{|\Gamma_t|}
  \left\langle \alpha(t), \beta(t_{(k)})\right\rangle. 
\end{equation}
By Lemma \ref{sum}, we have 
\begin{eqnarray*}
 \langle \alpha, d\beta \rangle_{L^2} 
 & = & \frac{(-1)^k}{k!}\sum_{s \in \orderedf{k-1}}
       \frac{1}{|\Gamma_s|}
       \sum_{t \in \ordered{k}_s}m(t)
       \left\langle \alpha(t), \beta(s)\right\rangle \\ 
 & = & \frac{1}{k!}\sum_{s \in \orderedf{k-1}} 
        \frac{m(s)}{|\Gamma_s|}\left\langle
        \left((-1)^{k} \sum_{t \in \ordered{k}_s} \frac{m(t)}{m(s)} 
         \alpha (t)\right), \beta(s)\right\rangle.
\end{eqnarray*}
This proves \eqref{adjoint_formula}.

Next we show that $d$ is bounded. 
For $\beta \in C_{\rho}^{k-1}(X)$, we have
\begin{eqnarray*}
  \left\Vert d\beta \right\Vert_{L^2}^2 
  & = & \frac{1}{(k+1)!} \sum_{t \in \orderedf{k}} \frac{m(t)}{|\Gamma_t|}
      \left\Vert d\beta(t) \right\Vert^2 \\
  & \leq & \frac{1}{k!} \sum_{t \in \orderedf{k}} \frac{m(t)}{|\Gamma_t|}
      \sum_{i=0}^{k}\left\Vert \beta(t_{(i)})\right\Vert^2. 
\end{eqnarray*}
By the argument that gives \eqref{rewrite} and the defining property
of $m$, the last expression becomes
\begin{equation*}
    \frac{k+1}{k!} \sum_{s \in \orderedf{k-1}} \frac{m(s)}{|\Gamma_s|}
       \left\Vert \beta(s)\right\Vert^2 \\
  = (k+1) \left\Vert \beta \right\Vert_{L^2}^2. 
\end{equation*}
Thus $d$ is bounded. Since $\delta$ is the adjoint of $d$,
$\delta$ is also bounded. This completes the proof. 
\end{proof}

The cohomology group of $X$ with coefficients in the representation 
$\rho$ is defined by 
\begin{equation*}
 H^k(X,\rho) = \left.
 \ker\left(d:C_{\rho}^k(X) \longrightarrow C_{\rho}^{k+1}(X)\right)
 \right/ \im\left(d:C_{\rho}^{k-1}(X) \longrightarrow 
 C_{\rho}^{k}(X)\right). 
\end{equation*}
If $X$ is contractible, then $H^k(X,\rho)\cong H^k(\Gamma,\rho)$, 
where $H^k(\Gamma, \rho)$ denotes the $k$-th cohomology group of
the group $\Gamma$ with coefficients in $\rho$.
Let $\Delta =d\delta + \delta d$.  
Then $\Delta$ is a bounded operator, and by elementary functional
analysis, we see that
\begin{eqnarray*}
 &&\ker\left(\Delta:C^k_{\rho}(X)\longrightarrow C^k_{\rho}(X)\right) \\
 &\cong&\ker\left(d:C^k_{\rho}(X)\longrightarrow C^{k+1}_{\rho}(X)\right)  
 \left/ \overline{
\im\left(d:C_{\rho}^{k-1}(X) \longrightarrow 
 C_{\rho}^{k}(X)\right)}\right..
\end{eqnarray*}
Therefore, if $\mathcal{H}$ is finite-dimensional, we have
\begin{equation}\label{dim_h_finite}
H^k(X,\rho) \cong \ker\Delta.
\end{equation}

Suppose that there is a constant $C>0$ such that 
$\langle \Delta \alpha, \alpha \rangle_{L^2} \geq C\Vert \alpha \Vert^2_{L^2}$
holds for all $\alpha \in C^k_{\rho}(X)$.
Then $\ker\Delta = 0$, and this implies $H^k(X,\rho)=0$ 
if $\mathcal{H}$ is finite-dimensional.
In fact, one can prove the vanishing of $H^k(X,\rho)$ in general
without referring to \eqref{dim_h_finite} as follows.
Since the spectrum bottom of $\Delta$ acting on $C^k_{\rho}(X)$ is positive, 
$\Delta$ has a bounded inverse $\Delta^{-1}$. 
Hence, any $\alpha \in C^k_{\rho}(X)$ decomposes as
 \begin{equation*}
  \alpha = \Delta \Delta^{-1} \alpha = d\delta \Delta^{-1} \alpha 
  + \delta d \Delta^{-1} \alpha,
 \end{equation*}
 which corresponds to the orthogonal decomposition 
 $C^k_{\rho}(X)= \ker d \oplus \overline{\im \delta}$.  
 Thus, if $\alpha$ is a cocycle, then we obtain 
$\alpha = d\delta \Delta^{-1} \alpha$, that is, $\alpha$ is a coboundary.
Therefore $H^k(X,\rho)=0$. 

We shall now write down a Bochner-Weitzenb\"ock-type formula for a $k$-cochain. 
This formula can be used to give a sufficient condition for the existence of 
a constant $C$ as above, 

\begin{Proposition}\label{linear_case}
For $\alpha \in C^k_{\rho}(X)$, $k\geq 1$, we have the following 
formula{\rm :}
\begin{equation}\label{linear_formula}
 \begin{split}
  & \Vert d\alpha \Vert_{L^2}^2 
     + \frac{k}{k+1}\Vert \delta \alpha \Vert_{L^2}^2  \\
& =  \frac{1}{k!}
     \sum_{s \in \orderedf{k-1}} \frac{1}{|\Gamma_s|} \left[
     \frac{1}{2} \sum_{(y,y') \in \overrightarrow{(\lk s)}(1)} 
     m(s,y,y') \Vert \alpha (s,y)-\alpha(s,y') \Vert^2 
     \right.\\
&  \phantom{\frac{1}{(k+1)!}} \left.
      - \frac{k}{k+1} \sum_{y \in (\lk s)(0)} m(s,y)
        \left\Vert \alpha (s,y) 
          - \sum_{y'\in (\lk s)(0)} \frac{m(s,y')}{m(s)}\alpha(s,y') 
    \right\Vert^2
  \right].
 \end{split}
\end{equation}
\end{Proposition}

\begin{proof}
By definition, we have
\begin{eqnarray*}
  \left\Vert \delta \alpha \right\Vert_{L^2}^2 
 & = &
   \frac{1}{k!}\sum_{s \in \orderedf{k-1}}\frac{m(s)}{|\Gamma_s|}
   \sum_{t,t' \in \ordered{k}_s} \frac{m(t)m(t')}{m(s)^2}
   \left\langle \alpha(t), \alpha(t') \right\rangle \\
 & = &
   \frac{1}{k!}\sum_{s \in \orderedf{k-1}}\frac{1}{|\Gamma_s|}
   \sum_{y,y' \in (\lk s)(0)} \frac{m(s,y)m(s,y')}{m(s)}
   \left\langle \alpha(s,y),\alpha(s,y')\right\rangle,
\end{eqnarray*}
and 
\begin{eqnarray*}
\lefteqn{\left\Vert d\alpha \right\Vert_{L^2}^2}\\ 
  &=&  \frac{1}{(k+2)!}\sum_{u \in \orderedf{k+1}}\frac{m(u)}{|\Gamma_u|}
  \left\langle \sum_{i=0}^{k+1}(-1)^i\alpha(u_{(i)})\ ,\ 
  \sum_{j=0}^{k+1}(-1)^j \alpha(u_{(j)}) \right\rangle \\
  &=&  \frac{1}{(k+2)!}\sum_{u \in \orderedf{k+1}}\frac{m(u)}{|\Gamma_u|}
  \left[ \sum_{i=0}^{k+1}\big\Vert\alpha(u_{(i)})\big\Vert^2 +
  2\sum_{i<j}(-1)^{i+j}
  \big\langle\alpha(u_{(i)}),\alpha(u_{(j)})\big\rangle\right]. 
\end{eqnarray*}
Note that for $i<j$, 
\begin{equation*}
\begin{split}
  & \alpha (u_{(i)}) = (-1)^{k-i}
  \alpha (\tau_{j,k+1}\circ \tau_{i,k}(u)_{(k)}), \\
  & \alpha (u_{(j)}) = (-1)^{k-j+1}
  \alpha (\tau_{j,k+1}\circ \tau_{i,k}(u)_{(k+1)})
\end{split}
\end{equation*}
hold. 
We now compute as in the proof of Proposition \ref{adjoint}:
\begin{equation}\label{dL2}
 \begin{split}
  \left\Vert d\alpha \right\Vert_{L^2}^2 
  & = \frac{1}{(k+2)!}\sum_{i=0}^{k+1} 
  \sum_{u' \in \tau_{i,k+1} (\orderedf{k+1})}\frac{m(u')}{|\Gamma_{u'}|}
  \left\Vert\alpha(u_{(k+1)}') \right\Vert^2 \\
  &\phantom{=} - \frac{2}{(k+2)!}\sum_{i<j}\ 
  \sum_{u' \in \tau_{j,k+1}\circ \tau_{i,k}(\orderedf{k+1})}
  \frac{m(u')}{|\Gamma_{u'}|}
  \left\langle \alpha(u_{(k)}'), \alpha(u_{(k+1)}')\right\rangle \\
  & = \frac{1}{(k+1)!} \sum_{s \in \orderedf{k-1}}\frac{1}{|\Gamma_s|}
  \sum_{t \in \ordered{k}_s} m(t) \left\Vert\alpha(t)\right\Vert^2 \\
  & \phantom{=} - \frac{1}{k!}\sum_{s \in \orderedf{k-1}}
  \frac{1}{|\Gamma_s|} \sum_{u \in \ordered{k+1}_s} m(u)
  \left\langle \alpha(u_{(k)}), \alpha(u_{(k+1)})\right\rangle \\
  & = \frac{1}{k!} \sum_{s \in \orderedf{k-1}}\frac{1}{|\Gamma_s|}
  \left[\frac{1}{k+1} \sum_{y \in (\lk s)(0)} m(s,y) 
  \left\Vert\alpha(s,y)\right\Vert^2 \right.\\
  & \phantom{=} - 
  \left. \sum_{(y,y') \in \overrightarrow{(\lk s)}(1)} m(s,y,y')
  \left\langle\alpha(s,y), \alpha(s,y')\right\rangle \right]. 
  \end{split}
\end{equation}
Noting that 
\begin{eqnarray*}
 &&  \frac{1}{k+1} \sum_{y \in (\lk s)(0)} m(s,y)
  \left\Vert\alpha(s,y)\right\Vert^2  - 
  \sum_{(y,y') \in \overrightarrow{(\lk s)}(1)} m(s,y,y')
  \left\langle\alpha(s,y), \alpha(s,y')\right\rangle \\
 &&  + \frac{k}{k+1}\sum_{y,y' \in (\lk s)(0)} 
    \frac{m(s,y)m(s,y')}{m(s)}
   \left\langle \alpha(s,y),\alpha(s,y')\right\rangle \\ 
& = &   \sum_{y \in (\lk s)(0)} m(s,y)
  \Vert \alpha(s,y) \Vert^2 
  - \sum_{(y,y') \in \overrightarrow{(\lk s)}(1)} m(s,y,y')
  \langle \alpha(s,y), \alpha(s,y')\rangle \\
 && - \frac{k}{k+1}
   \sum_{y \in (\lk s)(0)}m(s,y)\Bigg[\Vert \alpha(s,y) \Vert^2 
   -2 \sum_{y' \in (\lk s)(0)}\frac{m(s,y')}{m(s)}
    \langle \alpha(s,y), \alpha(s,y')\rangle \\
 && + \sum_{y',y'' \in (\lk s)(0)}\frac{m(s,y')m(s,y'')}{m(s)^2} 
    \langle \alpha(s,y'), \alpha(s,y'') \Bigg], 
\end{eqnarray*}
we obtain \eqref{linear_formula}.
\end{proof}

The expression in $[\quad]$ of \eqref{linear_formula} is
related to the first nonzero eigenvalue of the Laplacian 
of a finite graph as follows. 

For $s\in X(k-1)$ and $l=0,1$, we denote the set of $\mathcal{H}$-valued 
$l$-cochains on $\lk s$ by $C^l(\lk s, \mathcal{H})$.
Using the admissible weight $m$ of $X$, we define inner products on 
$C^l(\lk s,\mathcal{H})$, $l=0,1$, by
 \begin{equation*}
 \begin{split}
  & \langle f, g \rangle_{L^2} = \sum_{y \in (\lk s)(0)}m(s,y)
  \langle f(y), g(y) \rangle, \quad
   f,g \in C^0(\lk s, \mathcal{H}); \\
  & \langle \xi, \eta \rangle_{L^2} = \frac{1}{2}
  \sum_{(y,y') \in \overrightarrow{(\lk s)}(1)}
  m(s,y,y') \langle \xi(y,y'), \eta(y,y') \rangle, 
  \quad  \xi, \eta \in C^1(\lk s, \mathcal{H}). 
 \end{split}
\end{equation*}
With respect to these inner products, the Laplacian 
$\Delta_{\lk s}: C^0(\lk s, \mathcal{H}) \longrightarrow C^0(\lk s, \mathcal{H})$
takes the form of
\begin{equation*}
 (\Delta_{\lk s} f)(y)
 = f(y)- \sum_{y'; (y,y')\in \overrightarrow{(\lk s)}(1)}
\frac{m(s,y,y')}{m(s,y)} f(y').
\end{equation*}
We denote the first nonzero eigenvalue of $\Delta_{\lk s}$ by $\mu_1(\lk s)$.
If $\lk s$ is connected, the zero eigenspace of $\Delta_{\lk s}$
consists precisely of constant $0$-cochains.
Thus the variational characterization of $\mu_1(\lk s)$ 
reads
\begin{equation}\label{Rayleigh}
\mu_1(\lk s) = \inf \frac{||df||_{L^2}^2}{
||f-\overline{f}||_{L^2}^2},
\end{equation}
where the infimum is taken over all nonconstant $f\in C^0(\lk s, \mathcal{H})$
and
$$
\overline{f} = \sum_{y\in (\lk s)(0)} \frac{m(s,y)}{m(s)} f(y).
$$
It is an easy matter to verify that $\mu_1(\lk s)$ does not depend
on the choice of $\mathcal{H}$.

Let $\rho$ be a unitary representation of $\Gamma$ and $\alpha\in C^k_\rho(X)$.
Then for each $s \in \ordered{k-1}$, 
$(\lk s)(0)\ni y \mapsto \alpha(s,y) \in \mathcal{H}$
is an $\mathcal{H}$-valued $0$-cochain on $\lk s$.
An immediate consequence of \eqref{Rayleigh} is the following estimate:
 \begin{equation}
  \label{Rayleigh_quotient}
  \begin{split}
  & \frac{1}{2} \sum_{(y,y') \in \overrightarrow{(\lk s)}(1)} 
      m(s,y,y')
      \left\Vert \alpha (s,y)-\alpha(s,y') \right\Vert^2 \\
  \geq & \mu_1(\lk s)
  \sum_{y \in (\lk s)(0)}m(s,y)
  \left\Vert \alpha(s,y) 
     - \sum_{y'\in (\lk s)(0)} \frac{m(s,y')}{m(s)}\alpha(s,y')
  \right\Vert^2. 
  \end{split}
 \end{equation}

Now let $k\geq 1$.
Suppose that $\lk s$ is connected and 
\begin{equation}\label{spectral_condition3}
\mu_1(\lk s) > k/(k+1)
\end{equation}
for all $s \in \mathcal{F}(k-1)$.
Then there exists a constant $C>0$ 
such that $\langle \Delta \alpha, \alpha \rangle_{L^2}
\geq C\Vert \alpha \Vert^2_{L^2}$ for all $\alpha \in C^k_{\rho}(X)$.
In fact, since $\mathcal{F}(k-1)$ is a finite set, there exists a constant 
 $C'> k/(k+1)$ such that $\mu_1(\lk s) \geq C'$ for all 
 $s\in \mathcal{F}(k-1)$. 
 By \eqref{linear_formula} and \eqref{Rayleigh_quotient}, we obtain
 \begin{eqnarray*}
\lefteqn{\Vert d\alpha \Vert^2_{L^2} + \frac{k}{k+1}\Vert \delta \alpha
  \Vert^2_{L^2}} \\
  & \geq & \left(C' - \frac{k}{k+1}\right)\frac{1}{k!}
     \sum_{s \in \orderedf{k-1}} \frac{1}{|\Gamma_s|}
    \sum_{y \in (\lk s)(0)} m(s,y) \\
  &&\phantom{=}\times     \left\Vert \alpha (s,y) 
         - \sum_{y'\in (\lk s)(0)} \frac{m(s,y')}{m(s)}\alpha(s,y') 
   \right\Vert^2 \\
  & = & \left(C' - \frac{k}{k+1}\right)\frac{1}{k!}
     \sum_{s \in \orderedf{k-1}} \frac{1}{|\Gamma_s|}
     \sum_{y \in (\lk s)(0)} 
    \Bigg[ m(s,y) \Vert \alpha(s,y) \Vert^2  \\
  && \phantom{=} - \sum_{y' \in (\lk s)(0)}\frac{m(s,y)m(s,y')}{m(s)}
    \langle \alpha(s,y), \alpha(s,y')\rangle \Bigg] \\
  & = & \left(C' - \frac{k}{k+1}\right)
      \left((k+1)\Vert \alpha \Vert^2_{L^2} 
      - \Vert \delta \alpha \Vert^2_{L^2}\right). 
 \end{eqnarray*}
Adding $(C' - k/(k+1)) \Vert \delta \alpha \Vert^2_{L^2}$ 
to the both extreme sides, we obtain the desired result.

In summary, the condition \eqref{spectral_condition3} implies the vanishing of 
$H^k(X,\rho)$ for any unitary representation $\rho$. 
This result is due to Ballmann and \'Swi\c{a}tkowski \cite[Theorem 2.5]{BalSwi}.
In fact, they proved the vanishing of $L^2$-cohomology without assuming 
that the action of $\Gamma$ is cofinite.

We can relax the condition \eqref{spectral_condition3} as 

\begin{Proposition} 
\label{zuk's_formulation}
Suppose that $\lk s$ is connected for all $s \in X(k-1)$, and that
\begin{equation}
\label{spectral_condition4}
 \sum_{s \subset t}\mu_1(\lk s) > k
\end{equation}
holds for all $t \in X(k)$, where the sum is taken over all unordered
$(k-1)$-simplices in $t$. Then $H^k(X,\rho)$ vanishes.  
\end{Proposition}

When $k=1$, the result is due to \.Zuk \cite[Theorem 1]{Zuk1}. 

\begin{proof}
By a slight modification of the above computation, we have
 \begin{eqnarray*}
\lefteqn{\Vert d\alpha \Vert^2_{L^2} + \frac{k}{k+1}\Vert \delta \alpha
  \Vert^2_{L^2}} \\
  & \geq & \frac{1}{k!} \sum_{s \in \orderedf{k-1}}
    \frac{1}{|\Gamma_s|} \left(\mu_1(\lk s) - \frac{k}{k+1}\right) 
     \sum_{y \in (\lk s)(0)} 
     m(s,y) \Vert \alpha(s,y) \Vert^2  \\
  && \phantom{=} 
    -\left(2 - \frac{k}{k+1}\right)\Vert \delta \alpha \Vert^2_{L^2}, 
 \end{eqnarray*}
 since $\mu_1(\lk s)\leq 2$. Therefore,  
\begin{equation*}
   \Vert d\alpha \Vert^2_{L^2} + 2\Vert \delta \alpha
   \Vert^2_{L^2}
   \geq 
   \frac{1}{k!} \sum_{u \in \orderedf{k}}
    \frac{1}{|\Gamma_{u}|}
    \left(\mu_1(\lk u_{(k)}) - \frac{k}{k+1}\right) 
     m(u) \Vert \alpha(u) \Vert^2. 
\end{equation*}
Let $\frak{S}_{k+1}$ be the group of permutations on $k+1$ letters. 
Then, for any $\sigma \in \frak{S}_{k+1}$, $\sigma(\orderedf{k})$ is
also a representative set for the action of $\Gamma$ on $\ordered{k}$. 
Since the right-hand side of the above inequality is independent of the
choice of $\orderedf{k}$, we may rewrite the right-hand side as
 \begin{eqnarray*}
 & &  \frac{1}{k!(k+1)!}\sum_{\sigma \in \frak{S}_{k+1}} 
    \sum_{u \in \sigma(\orderedf{k})}
    \frac{1}{|\Gamma_{u}|}
    \left(\mu_1(\lk u_{(k)}) - \frac{k}{k+1}\right) 
     m(u) \Vert \alpha(u) \Vert^2 \\
 &=&  \frac{1}{k!(k+1)!}\sum_{\sigma \in \frak{S}_{k+1}} 
    \sum_{u \in \orderedf{k}}
    \frac{1}{|\Gamma_{\sigma(u)}|}
    \left(\mu_1(\lk \sigma(u)_{(k)}) - \frac{k}{k+1}\right) 
     m(\sigma(u)) \Vert \alpha(\sigma(u)) \Vert^2. 
 \end{eqnarray*} 
Note that $|\Gamma_{\sigma(u)}|$, $m(\sigma(u))$, and 
$\Vert \alpha(\sigma(u))\Vert$ are independent of 
$\sigma \in \frak{S}_{k+1}$. Changing the order of sum, we see
 that the last expression becomes 
\begin{equation*}
   \frac{1}{(k+1)!} \sum_{u \in \orderedf{k}} \frac{1}{|\Gamma_{u}|} 
    \sum_{i=0}^k \left(\mu_1(\lk u_{(i)})-\frac{k}{k+1}\right)
    m(u) \Vert \alpha(u) \Vert^2. 
\end{equation*}
By \eqref{spectral_condition4} and the finiteness of the set
 $\mathcal{F}(k)$, there exists a constant $C>0$ such that 
\begin{equation*}
 \Vert d\alpha \Vert_{L^2}^2 + 2 \Vert \delta \alpha \Vert_{L^2}^2
 \geq C \Vert \alpha \Vert_{L^2}^2
\end{equation*}
for all $\alpha \in C^k_{\rho}(X)$.  This completes the proof. 
\end{proof}

\end{document}